\documentclass[11pt]{amsart}
\usepackage{amssymb,mathrsfs,graphicx,enumerate}
\usepackage{amsmath,amsfonts,amssymb,amscd,amsthm,bbm}
\usepackage{kotex}
\usepackage{caption, subcaption}
\usepackage{comment}
\usepackage{graphicx,colortbl}
\usepackage{mathtools}
\usepackage{multirow}

\usepackage{xcolor}

\title[Emergent dynamics of the generalized Winfree model]{Emerging asymptotic patterns in a Winfree ensemble with higher-order couplings}

\author{Dongnam Ko}
\address{Department of Mathematics, The Catholic University of Korea, \\
Jibongro 43, Bucheon, Gyeonggido 14662, Republic of Korea\\
dongnamko@catholic.ac.kr}

\author{Seung-Yeal Ha}
\address{Department of Mathematical Sciences and Research Institute of Mathematics, \\ Seoul National University, Seoul 08826 \\
syha@snu.ac.kr}

\author{Jaeyoung Yoon}
\address{Department of Mathematical Sciences, \\ Seoul National University, Seoul 08826 \\
jyoung924@snu.ac.kr}

\newtheorem{theorem}{Theorem}[section]
\newtheorem{lemma}{Lemma}[section]
\newtheorem{corollary}{Corollary}[section]

\newtheorem{remark}{Remark}[section]

\newtheorem{definition}{Definition}[section]

\newcommand{\bbr}{\mathbb R}

\begin{document}
%%%%%%%%%%%%%%%%
\thanks{* Corresponding author.}
\date{\today}

\subjclass[2010]{} 

%\thanks{\textbf{Acknowledgment.} }

\thanks{\textbf{Acknowledgment.} The work of D. Ko was supported by the Catholic University of Korea, Research Fund, 2022, and by National Research Foundation of Korea (NRF-2021R1G1A1008559) and the work of S.-Y. Ha was supported by National Research Foundation of Korea (NRF-2020R1A2C3A01003881). Also, J. Yoon is grateful to the DFG-NRF International Research Training Group IRTG 2235 supporting the Bielefeld-Seoul graduate exchange programme.}

\begin{abstract}
The Winfree model is a phase-coupled synchronization model which simplifies pulse-coupled models such as the Peskin model on pacemaker cells. It is well-known that the Winfree ensemble with the first-order coupling exhibits discrete asymptotic patterns such as incoherence, locking and death depending on the coupling strength and variance of natural frequencies. In this paper, we further study higher-order couplings which makes the dynamics more close to the behaviors of the Peskin model. For this, we propose several sufficient frameworks for asymptotic patterns compared to the first-order coupling model. Our proposed conditions on the coupling strength, natural frequencies and initial data are independent of the number of oscillators so that they can be applied to the corresponding mean-field model. We also provide several numerical simulations and compare them with analytical results. 
\end{abstract}
\maketitle \centerline{\date}

%\tableofcontents
\section{Introduction} \label{sec:1}
\setcounter{equation}{0}
Synchronization denotes the adjustment of rhythms in an oscillatory complex system. After novel approaches by Arthur Winfree and Yoshiki Kuramoto in almost half century ago, it has been extensively studied in diverse scientific disciplines such as applied mathematics, biology, control theory and statistical physics, etc. 

In literature, there are two types of synchronization models for coupled oscillators. 
The first type is called {\it``phase-coupled model"}. This model assumes that the amplitude variations of oscillator's states are assumed to be negligible, and it focuses only on the variations of phases. The models introduced by Winfree and Kuramoto belong to this category. 
In contrast, there are pulse-coupled models such as an integrate-and-fire model, Peskin model for pacemaker cells, etc. This latter type of models might be able to explain synchronization in real neurons, however, it is very difficult to perform rigorous mathematical analysis. This is why phase-coupled models are more frequently employed in the study of synchronization. The Peskin model \cite{M-S,Peskin} requires discrete jumps in dynamics, while its asymptotic properties are recently discovered in some special cases \cite{Akhmet,Z-D-G-B-Z}. Another famous example is a network of FitzHugh–Nagumo oscillators, where its synchronization behavior has been extensively studied by researchers from mathematics, biology and physics \cite{A-B-F-H-K, N-H,P-F,S-T-N}. 

Analytical study on synchronization has been intensively complemented on the studies of oscillators in $\mathbb T^N$ with smooth interactions \cite{A-B,D-B1,Ku,P-R-K,St} to avoid chaotic dynamical patterns. 
For a positive coupling strength, the firing phenomenon of one oscillator affects other oscillator's phase. However, when the coupling strength is sufficiently small, the dynamics of each phase is basically governed by the corresponding natural frequency in the presence of coupling (called emergence of incoherence). On the opposite extreme regime in which the coupling strength is sufficiently large, the ensemble of phases approaches to an equilibrium and oscillators will not rotate any more. We call this emergent pattern as {\it death}. The remaining case is the intermediate regime in which the coupling strength is not that small and not that large so that all oscillators have the same rotation number. This pattern is called {\it locking}.

To put our discussion in a proper setting, we begin with the description of the Winfree model. Let $\theta_i = \theta_i(t)$ be the phase of the $i$-th Winfree oscillator. In a pulse-coupled setting, it is common to consider that one individual `fires', when its phase reaches a certain value, namely, zero. Therefore, the influence of one oscillator to others should be accumulated near the zero value. To model such dynamics, we use a smooth but sharp bump-like function which we call it as ``influence function" $I_n$ with an order $n$. Among phase-coupled models, the Winfree model \cite{Winfree} is, from the beginning, built to approximate pulse-coupled synchronization using a smoothly approximated influence function $I_n$ for the constant multiple of Dirac-delta function:
\begin{equation} \label{A-0}
2\pi \delta(\theta) \approx a_n (1 + \cos \theta)^{n} =: I_n, \quad n \gg 1, \quad |\theta| \ll 1.
\end{equation}
Here $n$ is the order of trigonometric coupling, which is a positive integer, while the coefficient $a_n$ is a positive normalizing constant chosen to satisfy
\[  \int_{-\pi}^{\pi}I_n(\theta)d\theta= 2\pi. \]
In this setting, the Winfree model with the influence function \eqref{A-0} reads as 
\begin{equation}\label{A-1}
\displaystyle \dot{\theta_i}=\nu_i+  \frac{\kappa}{N} \sum_{j=1}^N   I_n(\theta_j) S(\theta_i), \quad i \in [N]:= \{1, \cdots, N \},
\end{equation}
where $\nu_i$ is the natural frequency of the $i$-th oscillator, and $(I_n, S)$ is the pair of influence-sensitivity functions of order $n \in {\mathbb Z}_+$:
\begin{equation} 
\begin{cases} \label{A-2}
\displaystyle I_n(\theta):= a_n (1+\cos \theta )^n = 2^n a_n \cos^{2n}\frac{\theta}{2}, \quad  S(\theta):=-\sin\theta, \vspace{0.5em}\\
\displaystyle  a_n =\frac{2n(2n-2)\cdots 2}{2^n (2n-1)(2n-3)\cdots 1} = \frac{(2n)!!}{2^n(2n-1)!!}.
\end{cases}
\end{equation}
\begin{figure}
		\centering
		\begin{subfigure}[ht]{0.4\textwidth}
\centering
			\includegraphics[width=1.15\textwidth]{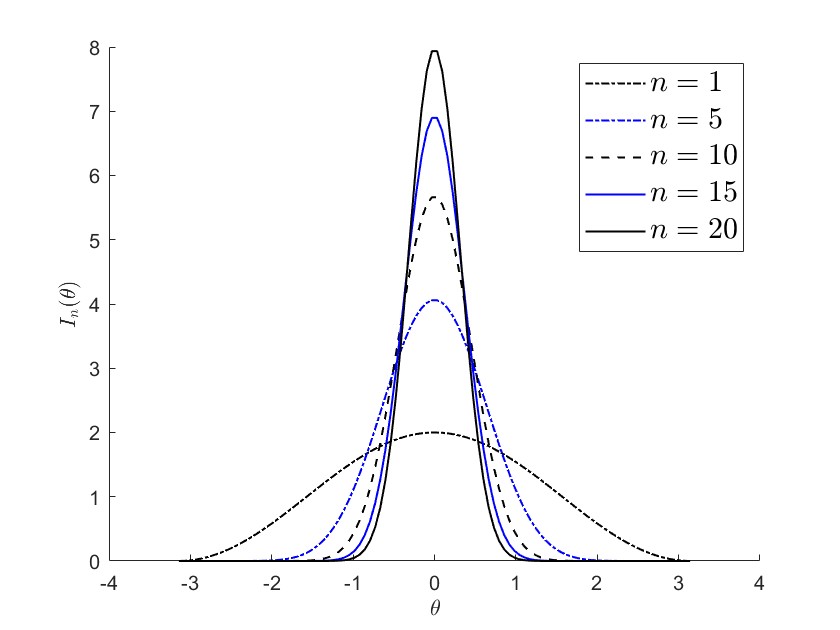}\caption{Graphs of $I_n$ with different $n$}
		\end{subfigure}
\quad
		\begin{subfigure}[ht]{0.4\textwidth}
\centering
			\includegraphics[width=1.15\textwidth]{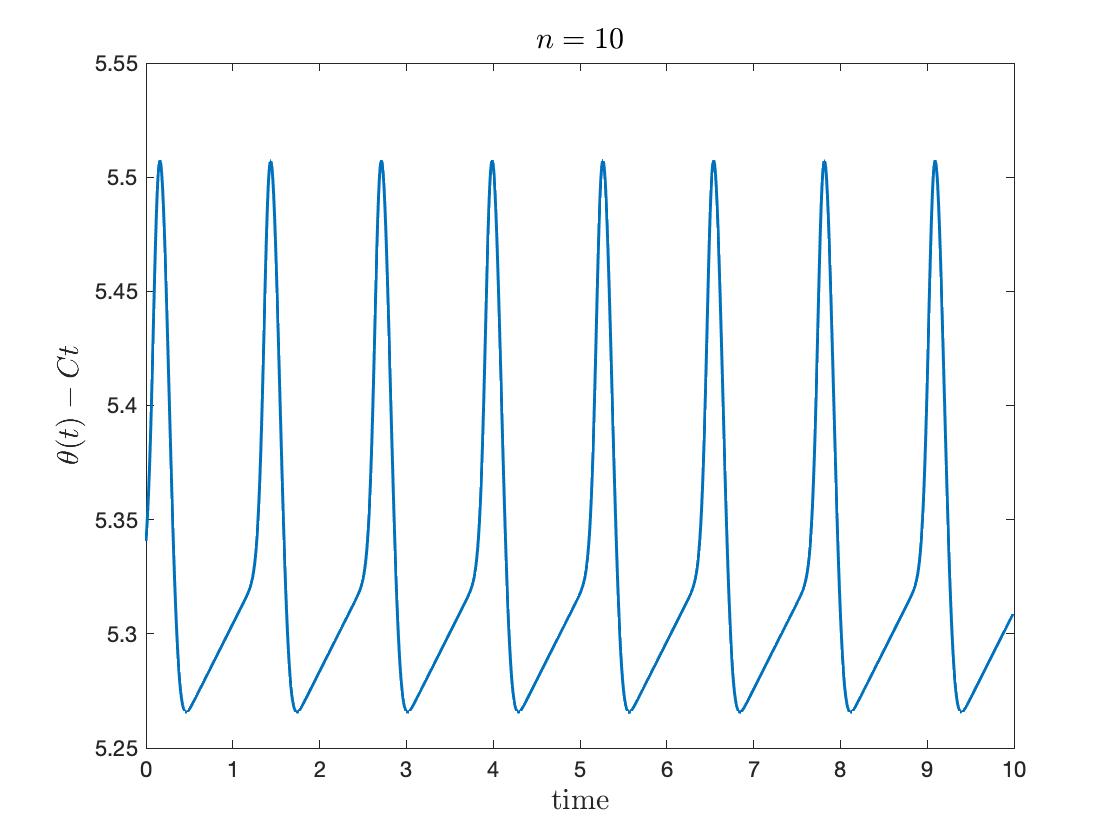}\caption{Adjusted phase dynamics for $n= 10$} 	
			\end{subfigure}
	\caption{Single-oscillator dynamics for $n=10$}
		\label{Fig-1}
	\end{figure}

From now on, we call \eqref{A-1}--\eqref{A-2} as {\it``the Winfree model with order $n$"}. As we can see in Figure \ref{Fig-1}(A), the maximal height of $I_n$ grows with the order of $\mathcal{O}(\sqrt{n})$ (check Lemma \ref{L3.1} for details), whereas the effective size of the support of $I_n$ decreases as $n$ increases. As $n \to \infty$, the influence function of order $n$ approximates Dirac-delta function so that the resulting dynamics might behave the same way as the pulse-coupled one. In Figure \ref{Fig-1}(B), we can see the temporal evolution of phase over time while we subtract a kind of group velocity to the phase graph in order to observe the pulse-like behavior significantly. 

In previous literature \cite{A-S, M-S}, the first-order pair $(I_1, S)$:
\[
 I_1(\theta) = 1+\cos\theta \quad \mbox{and} \quad S(\theta) = -\sin\theta
\]
has been extensively used in the mathematical analysis on the emergent dynamics of a Winfree emsemble. 
To introduce definitions on asymptotic patterns, here we adopt a classical concept, namely {\it rotation number}, which is commonly used in literature of phase-coupled dynamics. Let $\Theta = \{ \theta_i \}_{i \in [N]}$ be a collection of phases. For each $i \in [N]$, we define the rotation number $\rho_i$ as an asymptotic frequency (\cite{K-H}):
\[ \rho_i := \lim_{t \to \infty} \frac{\theta_i(t)}{t}. \]
if the right-hand side exists. Equipped with rotation numbers, we can classify asymptotic patterns in the Winfree model as follows.
\begin{definition}\label{D1.1} \cite{A-S}
Let $\Theta = \Theta(t)$ be a time-dependent oscillator ensemble.
\begin{enumerate}
\item The ensemble $\Theta(t)$ exhibits \emph{death} asymptotically if all the rotation numbers are zero:
\[ \rho_i = 0, \quad \forall~ i \in [N]. \]
\item The ensemble $\Theta(t)$ exhibits \emph{locking} asymptotically if all the rotation numbers are equal to a nonzero constant:
\[ \rho_i = \rho \neq 0, \quad \forall~i \in [N],  \]
\item The ensemble $\Theta(t)$ exhibits \emph{incoherence} if all the rotation numbers for oscillators with different natural frequency are different:
\[ \rho_i \neq \rho_j, \quad \forall~(i,j) \in [N] \times [N]~\mbox{such that}~\nu_i \neq \nu_j. \]
\end{enumerate}
\end{definition}
%As observed in \cite{H-P-R}, the Winfree model \eqref{A-1} is a gradient flow only for $n = 1$ while one cannot use gradient structure with $n > 1$ any more. This mathematical difficulty, fortunately, can be avoided using the machinery employed in \cite{H-K-P-R, H-P-R}. 
\vspace{0.2cm}
The purpose of this paper is to investigate the emergent dynamics of \eqref{A-1}--\eqref{A-2} depending on the order $n$. 
The bifurcation simulation in \cite{A-S} of $S$ and $I_n$ with $n=10$ presents a similar diagram qualitatively to the one with $n=1$. However, their critical values of $\kappa$ for incoherence becomes much smaller. This phenomenon may indicate that the bifurcation of collective dynamics in the original Winfree model \eqref{A-1} may not be preserved in the limit $n \to \infty$ since the critical coupling strength for bifurcation could tend to $0$. This generates the following question: \newline
\begin{quote}
``What is a proper normalization factor for the interaction kernels $I_n$ and $S$ which we can still see the bifurcation diagram between incoherence, locking and death, as $n \to \infty$?"
 \end{quote}
 \vspace{0.2cm}
 
 Next, we briefly discuss our main results. First, we present a sufficient framework for incoherence in terms of natural frequency and coupling strength:
 \[  \nu_i \neq \nu_j, \quad  0 \leq \kappa <  \kappa^{ij}_{\mathrm{inc}} := \frac{ |\nu_i - \nu_j|}{2^n a_n} \sim {\mathcal O}\Big( \frac{1}{\sqrt{n}}\Big). \]
 Under these conditions, there exists a positive constant $\omega_{ij}^{\infty}$ for each pair of $i$ and $j$ such that 
\[ \inf_{0 \leq t < \infty} |\dot{\theta}_i(t)-\dot{\theta}_j(t)| \geq \omega_{ij}^{\infty} > 0. \]
Note that $\kappa^{ij}_{\mathrm{inc}} $ does not depend on $N$. We refer to Theorem \ref{T3.1} and Remark \ref{R3.2} for details. 

Second, we provide a sufficient framework leading to {\color{black}death.} Let $\beta_n$ be the minimum point of $SI_n$ in $(0, \pi)$ determined later by the relation \eqref{D-4-1}, and we assume that parameters $\alpha_n, \kappa$ and initial data satisfy
\[
\beta_n < \alpha_n < \pi,  \quad  |\theta_i^0| < \alpha_n, \quad  \kappa>   \frac{ \max_{i} |\nu_i|}{|(SI_n)(\alpha_n)|} = \kappa_{d,n}(\alpha_n). 
\]
Then, all the rotation numbers are zero (see Theorem \ref{T4.1}):
\[ \rho_i = 0, \quad i \in [N]. \] 
Here the critical coupling strength $\kappa_{d,n}(\beta_n)$, where $\beta_n$ is the maximizer of $|SI_n|$, is the order of ${\mathcal O}(1)$. This is consistent with the phase diagram in \cite{A-S} that $\kappa_{d,n}(\beta_n)$ is not affected much by $n$. 

Third, we deal with a sufficient condition for locking in a homogeneous Winfree ensemble, where the whole oscillators have the same natural frequency $\nu > 0$. For some restricted initial phase configuration, if coupling strength is sufficiently small 
\[ 0 < \kappa < \frac{\nu}{2^{n+1} a_n} \sim {\mathcal O}\Big(\frac{\nu}{\sqrt{n}}\Big), \]
then all the rotation numbers are equal to $\nu$:
\[  \rho_i = \nu, \quad i \in [N]. \]
We refer to Theorem \ref{T5.1} for details. When $\kappa$ becomes large, it tends to death by Theorem \ref{T4.1}. \newline
 
The rest of this paper is organized as follows. In Section \ref{sec:2}, we briefly review basic structures for the interaction pair $(I_n, S)$ and recall previous results on the emergent dynamics of the Winfree model with order one. In Section \ref{sec:3} and Section \ref{sec:4}, we provide a rigorous analysis for the emergence of incoherence and death, respectively. In Section \ref{sec:5}, we present sufficient frameworks leading to complete and partial phase-lockings for a homogeneous Winfree emsemble. In Section \ref{sec:6}, we provide several numerical examples and compare them with analytical results in proceeding sections. Finally, Section \ref{sec:7} is devoted to a brief summary of our main results and discussions on some remaining issues. 

\hspace{0.5cm}

\noindent {\bf Notation}:~For simplicity, we employ the following handy notation:
\begin{align*}
\begin{aligned}
& \Theta :=  \{ \theta_i \}, \quad  \mathcal{V}: =  \{ \nu_i \}, \quad  \|\Theta \|_{p} := \Big( \sum_{i=1}^{N} |\theta_i|^p \Big)^{\frac{1}{p}}, \quad 1 \leq p < \infty \quad  \|\Theta \|_{\infty} := \max_{1\leq i \leq N} |\theta_i|, \\
& {\mathcal D}(\Theta) := \max_{1\leq i,j\leq N}|\theta_i -\theta_j|,\quad {\mathcal D}(\mathcal{V}) := \max_{1\leq i,j\leq N}|\nu_i - \nu_j|.
\end{aligned}
\end{align*}
For a function $F: \bbr_+ \to \bbr$, we denote the $L^p(\bbr_+)$-norm of $F$ by 
\[ \| F \|_{L^p(\bbr_+)} := \begin{cases}
\Big( \int_0^{\infty} |F(t)|^p dt  \Big)^{\frac{1}{p}}, \quad & 1 \leq p < \infty, \\
\sup_{0 \leq t < \infty} |F(t)|, \quad & p = \infty.
\end{cases}
\]
\section{Preliminaries}\label{sec:2}
\setcounter{equation}{0}
In this section, we discuss the reinterpretation of \eqref{A-1} as a generalized Adler equation, and we recall some structural conditions for influence-sensitivity function pair $(I_n, S)$. We also present previous results on the emergent dynamics of the Winfree model \eqref{A-1} with order $1$. 
\subsection{The Winfree model with order $n$} \label{sec:2.1} In this subsection, we briefly discuss how the Winfree model with order $n$ and the Kuramoto model can be viewed as generalized Adler equation.  As noticed in \cite{C-H-J-K}, the Adler equation:
\begin{equation} \label{B-0}
{\dot \theta} = \nu - \kappa \sin \theta, \quad t > 0,
\end{equation}
 plays a key role in synchronization dynamics. More precisely, we consider the Kuramoto model for two-oscillators:
\[
{\dot \theta}_1 = \nu_1 + \frac{\kappa}{2} \sin(\theta_2 - \theta_1), \quad  {\dot \theta}_2 = \nu_2 + \frac{\kappa}{2} \sin(\theta_1 - \theta_2).
\]
In this case, the phase difference $\theta = \theta_1 - \theta_2$ satisfies \eqref{B-0}. On the other hand, the Kuramoto model with $N \geq 3$:
\[
{\dot \theta}_i = \nu_i + \frac{\kappa}{N} \sum_{j=1}^{N} \sin(\theta_j  - \theta_i)
\]
can also be written as a generalized Adler equation with state dependent coupling strength $\kappa = \kappa(\Theta)$ (see \cite{C-H-J-K}) when the geometric shape of a phase configuration is confined in a half circle:
\[
\dot{\theta}_{ij} = \nu_{ij} - \kappa \Big(  \sum_{l \not = i, j}  C^l_{ij}(\Theta) \Big) \sin \theta_{ij}, 
\]
where differences $\theta_{ij},~\nu_{ij}$ and mean-field factor $C^l_{ij}$ are given as follows.
\[
\theta_{ij}:= \theta_i - \theta_j, \quad \nu_{ij}:= \nu_i - \nu_j, \quad C^l_{ij}(\Theta):= 1 - \frac{\cos (\frac{\theta_{li}}{2} + \frac{\theta_{lj}}{2})}{\cos (\frac{\theta_{ji}}{2}  )}. 
\]
Similarly,  \eqref{A-1}  can cast as a Adler type equation. To see this, we set $I_{n,c}(\Theta)$ as the average of $\{ I(\theta_i) \}$:
\[
I_{n,c}(\Theta): = \frac{1}{N} \sum_{k=1}^{N} I_n(\theta_k). 
\]
Then, the Winfree model \eqref{A-1}--\eqref{A-2} is rewritten as follows:
\[
\dot{\theta_i}=\nu_i - \kappa I_{n,c}(\Theta) \sin(\theta_i), \quad i \in [N].
\]
In this way, the Kuramoto model and the Winfree model with order $n$ are both the special cases for a generalized Adler equation, and the only difference lies in the functional form of coupling strength on $\Theta$:
\[ 
\mbox{Kuramoto oscillator}:~\kappa \sum_{l \not = i, j}  C^l_{ij}(\Theta), \qquad \mbox{Winfree oscillator}:~ \frac{\kappa}{N} \sum_{k=1}^{N} I(\theta_k), 
\]
which generates notable differences in asymptotic dynamics of both models. 

\subsection{Structural conditions on $(I_n, S)$} \label{sec:2.2}
Next, we list sufficient conditions for the interaction kernels $(I_n, S)$. Motivated by the specific example \eqref{A-2}, the following conditions were proposed in \cite{H-P-R}: 
\begin{itemize}
\item ($\mathcal{A}1$)~(Periodicity and parity conditions):  for $\theta \in \mathbb{R},$
\begin{equation}\label{B-1}
S(\theta+2\pi)=S(\theta),\quad S(-\theta)=-S(\theta),\quad I_n(\theta+2\pi)=I_n(\theta),\quad I_n(-\theta)=I_n(\theta).
\end{equation}
\item ($\mathcal{A}2$)~(Geometric shapes): there exist positive constants $\theta_{*}$ and $\theta^{*}$ satisfying 
\[ 0<\theta_{*}<\theta^{*} \leq \pi\]
such that
\begin{equation} 
\begin{cases} \label{B-2}
\displaystyle S \leq 0 \quad\text{on }~[0,\theta^*] \quad\text{and}\quad S' \leq 0, \quad S'' \geq 0 \quad\text{on }~ [0,\theta_*],\\
\displaystyle I_n \geq 0,\quad I_n'\leq 0 \quad\text{on }~[0,\theta^*] \quad\text{and}\quad I_n'' \leq 0 \quad\text{on }~ [0,\theta_*],\\
\displaystyle (SI_n)' < 0 \quad\text{on }~(0,\theta_*) \quad\text{and}\quad (SI_n)'>0 \quad\text{on }~(\theta_*,\theta^*).
\end{cases}
\end{equation}
\end{itemize}
\begin{remark} Our choice \eqref{A-2} satisfies ($\mathcal{A}1$)--($\mathcal{A}2$) with
\[ \theta_* = \cos^{-1} \Big( \frac{n}{n +1} \Big), \quad \theta^* = \pi. \]
\end{remark}

\subsection{Previous results} \label{sec:2.3}
In this subsection, we briefly review previous results in \cite{H-P-R} on the emerging patterns for the Winfree model with order $1$. First, we recall a sufficient framework leading to the complete oscillator death in which all the rotation numbers are zero. 
\begin{theorem}\label{T2.1}
\emph{ (Emergence of death \cite{H-P-R})} 
Suppose that parameters and initial data satisfy 
\begin{equation} \label{B-3}
0 < \alpha < \pi, \qquad \kappa > -\frac{\| {\mathcal V} \|_{\infty}}{S(\alpha)I_1(\alpha)}, \qquad  \max_{i \in [N]} |\theta_i^0|\leq\alpha,
\end{equation}
and let $\Theta = \Theta(t)$ be a solution to \eqref{A-1}--\eqref{A-2}. Then, the complete oscillator death emerges:
\[  \rho_i = 0, \quad i \in [N].  \]
\end{theorem}
\begin{remark} Under the setting \eqref{B-3}, one can show that there exists an equilibrium state $\Theta^{\infty}$ such that 
\[ \lim_{t \to \infty} \| \Theta(t) - \Theta^{\infty} \|_{\infty} = 0. \]
Clearly, this is much stronger than the fact that the ensemble exhibits death. 
\end{remark}
In next theorem, we consider complete phase-locking. %For $\alpha \in (0, \pi)$, we set $\alpha^{\infty}$
\begin{theorem}
\emph{(Emergence of locking \cite{H-K-P-R})} \label{T2.2}
The following assertions hold.
\begin{enumerate}
\item
(Identical ensemble):~If parameters and initial data satisfy
\[ 
\begin{cases}
\displaystyle \nu_i = \nu, \quad \forall~i \in [N], \quad 0 < \kappa < \frac{\pi}{8} \nu, \quad  0 < \alpha <  {\color{black}\Big(\frac{\pi}{4} - \frac{2\kappa}{\nu} \Big)}, \\
\displaystyle {\mathcal D}(\Theta^0)\leq \alpha  \exp \Big[  -\frac{\kappa}{\nu - 2 \kappa} \Big( 1 + \alpha - \frac{\pi}{2} \Big)      \Big ], 
\end{cases}
\]
and let $\Theta(t)$ be a solution to \eqref{A-1}-\eqref{A-2}.  Then, there exist positive constants $\beta_i$ and $\Lambda_i$ such that  
\[ \beta_1  {\mathcal D}(\Theta^0) e^{-\Lambda_1 t} \leq {\mathcal D}(\Theta(t)) \leq \beta_2  {\mathcal D}(\Theta^0) e^{-\Lambda_2 t} \leq \alpha,\quad t \geq 0. \] 
\item
(Non-Identical ensemble):~If parameters and initial data satisfy
\begin{align*}
&\nu > 0, \quad \kappa > 0, \quad {\mathcal D}(\Theta^0) < \alpha, \quad {\color{black} 0 < \alpha <  \Big(\frac{\pi}{8} - \frac{\kappa}{\nu} \Big)},\\
& {\color{black}0 < \delta_{\nu} < \nu-6\kappa}, \quad \nu - \delta_{\nu} \leq \nu_i \leq \nu + \delta_{\nu}, \quad \forall~ i \in [N].
\end{align*} 
Then, for a solution $\Theta(t)$ to \eqref{A-1}--\eqref{A-2}, one has 
\[ \sup_{0 \leq t < \infty} {\mathcal D}(\Theta(t))\leq 3\alpha. \]
\end{enumerate}
\end{theorem}
\begin{remark}
\begin{enumerate}
\item
The explicit functional relations for $\beta_i$ and $\Lambda_i$ and partial phase-locking can be found in \cite{H-K-P-R}. 
\item
We also refer to other related results \cite{O-T-K, O-K-T} for the Winfree model.
\item
Phase-locking of the Kuramoto oscillators has been extensively studied in literature, to name a few, we refer to \cite{A-R, B-C-M, B-D-P, C-S, D-X, D-B1, D-B2, H-K-R, H-K-P-Z, H-R, J-M-B, V-M0, V-M1, V-M2}.
\end{enumerate}

\end{remark}
In the following three sections, we prove the emergence of asymptotic patterns discussed in Definition \ref{D1.1} with higher-order couplings, which correspond to Theorem \ref{T2.1} and Theorem \ref{T2.2}.
\section{Emergence of incoherence}\label{sec:3}
\setcounter{equation}{0}
In this section, we present a sufficient framework leading to incoherent state. For this, we show that there exists a pair of indices $(i_0, j_0) \in [N] \times [N]$ such that frequency difference $|{\dot \theta}_{i_0} - {\dot \theta}_{j_0}|$ has a positive lower bound, which shows that these two oscillators cannot be entrained. First, we estimate the bounds of interaction functions in the following lemma.
\begin{lemma}\label{L3.1}
The following assertions hold.
\begin{enumerate}
\item
The normalized constant $a_n$ in \eqref{A-2} is the order of $\sqrt{n}$ asymptotically for $n \gg 1$:
\[ a_n \frac{2^n}{\sqrt{n \pi}} \searrow 1 \qquad\text{as}\quad n \to \infty. \]
\item
Interaction function pair $(I_n, S) $ satisfies
\begin{align*}
\begin{aligned}
& (a)~ \|I_n\|_{L^{\infty}(\bbr)} = 2^n a_n \sim \mathcal O( \sqrt{n} ), \quad  \|S\|_{L^{\infty}(\bbr)} = 1, \\
&(b)~ \|SI_n\|_{L^\infty(\bbr)} =  a_n\frac{\sqrt{2n+1}}{n+1}\left(\frac{2n+1}{n+1}\right)^n\sim \mathcal O\left(1\right), \\
&(c)~ \| I^{\prime} \|_{L^\infty(\bbr)}  =  a_n  \Big( \frac{2n-1}{n}  \Big)^{n-1} \sqrt{2n - 1} \sim{\mathcal O}\left(n\right),
\end{aligned}
\end{align*}
as $n \to \infty$. 
\end{enumerate}
\end{lemma}
\begin{proof}  We set 
\[  b_n := \frac{a_n 2^n}{\sqrt{n \pi}}. \]

\noindent(i)~First, we show that 
\[  \lim_{n \to \infty} b_n  = 1. \]
For this, we use \eqref{A-2} and  Stirling's formula:
\begin{equation} \label{C-2}
2^n a_n = \frac{(2n)!!}{(2n-1)!!} = \frac{((2n)!!)^2}{(2n)!} = \frac{2^{2n}((n)!)^2}{(2n)!} \quad \mbox{and} \quad  \lim_{n \to \infty} \frac{\sqrt{2\pi n}}{n!}\left(\frac{n}{e}\right)^n = 1 
\end{equation}
to get
\begin{align*}
\lim_{n \to \infty} \frac{1}{\sqrt{n\pi}}\frac{(2n)!!}{(2n-1)!!} = \lim_{n \to \infty} \frac{1}{\sqrt{n\pi}}\frac{2^{2n}((n)!)^2}{(2n)!} = \lim_{n \to \infty} \frac{1}{\sqrt{n\pi}}\frac{2^{2n}\cdot 2\pi n \cdot n^{2n}e^{-2n}}{\sqrt{4\pi n}\cdot  2^{2n}n^{2n}e^{-2n}} = 1.
\end{align*}
Therefore, we have
\[ \lim_{n \to \infty} \frac{2^n a_n}{\sqrt{n \pi}} = 1. \]
Second, we show that the sequence $(b_n)$ is strictly decreasing: By \eqref{C-2}, we have
\[
\frac{b_{n+1}}{b_n}=\frac{(2n+2)!!}{(2n+1)!!}\frac{(2n-1)!!}{(2n)!!}\frac{\sqrt{n}}{\sqrt{n+1}} = \frac{2n+2}{2n+1}\frac{\sqrt{n}}{\sqrt{n+1}} = \frac{\sqrt{n+1}\sqrt{n}}{n+1/2}<1, \quad \forall~n \geq 1. 
\]

\vspace{0.2cm}

\noindent(ii)~Recall that 
\begin{equation} \label{C-2-0}
I_n(\theta) = a_n (1 + \cos \theta)^{n} = 2^n a_n \left(\cos \frac{\theta}{2}  \right)^{2n} \quad\text{and}\quad S(\theta) =-\sin\theta. 
\end{equation}

(a)~The bounds of $I_n$ and $S$ are obvious from \eqref{C-2-0} and the estimate of $a_n$. \newline

(b)~We set a $2\pi$-periodic function $f$ as 
\[ f(\theta) = I_n(\theta) S(\theta) = -a_n (1 + \cos \theta)^{n} \sin \theta. \]
By direct calculation, we have
\[ f^{\prime}(\theta) = a_n (1 + \cos \theta)^n \Big[ n - (n+1) \cos \theta \Big]. \]
We set $\theta_*$ to be the minimizer of $f$ in $[-\pi,\pi]$, i.e., 
\[  \cos \theta_* = \frac{n}{n + 1} \quad \text{and} \quad \sin \theta_* = \sqrt{  1 - \Big(  \frac{n}{n+1}  \Big)^2 } = \frac{\sqrt{2n + 1}}{n + 1}.   \]
It follows from the graph of $|f|$ that it has a maximum at $\theta_*$ and $-\theta_*$ in the domain of $[-\pi,\pi]$. Its maximum value is 
\[ \|SI_n \|_{L^{\infty}(\bbr)} = a_n \Big(  1 + \cos \theta_* \Big)^{n} \sin \theta_* = a_n \Big( 1  + \frac{n}{n + 1} \Big)^n \frac{\sqrt{2n + 1}}{n + 1} = 
a_n \Big( \frac{2n + 1}{n + 1} \Big)^n \frac{\sqrt{2n + 1}}{n + 1} . \]
\noindent (c)~We will use the same argument as in (b). More precisely, we set 
\[ g(\theta) = I_n^{\prime}(\theta) =  -a_n n (1 + \cos \theta)^{n-1} \sin \theta. \]
By direct calculation,
\[
g^{\prime}(\theta) =  a_n n ( 1 + \cos \theta)^{n-1} \Big[  (n-1) - n \cos \theta  \Big ].
\]
Thus, $|g|$ has a maximum at the value ${\tilde \theta}_*$ such that 
\[ \cos {\tilde \theta}_* = \frac{n-1}{n}, \quad \sin {\tilde \theta}_*  =  \frac{\sqrt{2n - 1}}{n},\]
i.e., 
\[ |g({\tilde \theta}_*)| = n a_n   \Big ( 1 +  \frac{n-1}{n} \Big)^{n-1}  \frac{\sqrt{2n - 1}}{n} =  n a_n  \Big( \frac{2n-1}{n}  \Big)^{n-1} \frac{\sqrt{2n - 1}}{n} = 
{\color{black}\mathcal{O}(n)}.  \]
\end{proof}
\begin{remark}\label{R3.1}
It follows from Lemma \ref{L3.1} that 
\[ 1\leq b_n \leq b_1 =   \frac{2}{\sqrt{\pi}}  a_1 = \frac{2}{\sqrt{\pi}}   \quad \forall~n \geq 1 \]
which is equivalent to 
\begin{equation} \label{C-2-1}
\frac{\sqrt{n \pi} }{2^n}\leq a_n \leq \frac{\sqrt{n}}{2^{n-1}},\quad n \geq 1.
\end{equation}
\end{remark}
\vspace{0.5cm}
Now, we present our first main result on the emergence of incoherent state. If the coupling strength vanishes ($\kappa =0$), oscillators disperse according to their natural frequencies. Similarly,  for a small coupling strength $\kappa \ll 1$, oscillators with different natural frequencies will disperse. We assume that all the natural frequencies are different and set 
\begin{align*}
\begin{aligned}
& \kappa^{ij}_{\mathrm{inc}} := \frac{|\nu_i - \nu_j|}{2^{n+1} a_n}, \quad i, j \in [N], \quad  \kappa_{\mathrm{inc}} := \frac{\min_{i \neq j}  |\nu_i - \nu_j|}{2^{n+1} a_n}, \\
& \omega_{ij}^{\infty} := |\nu_i - \nu_j| -\kappa 2^{n+1} a_n, \quad \omega_m^{\infty} := \min_{i \neq j} \omega_{ij}^{\infty}. 
\end{aligned}
\end{align*}
\begin{theorem}\label{T3.1}
\emph{(Emergence of incoherence)}
Suppose that system parameters satisfy 
\begin{equation} \label{C-3}
\nu_i > \nu_j, \quad 0 \leq \kappa <  \kappa_{\mathrm{inc}}
\end{equation}
and let $\Theta$ be a global solution to \eqref{A-1}--\eqref{A-2}.  Then, we have
\[ \inf_{0 \leq t < \infty} |\dot{\theta}_i(t)-\dot{\theta}_j(t)| \geq \omega_{ij}^{\infty} > 0. \]
\end{theorem}
\begin{proof}
Recall a functional form $(I_n, S)$:
\begin{equation} \label{C-3-0-1}
I_n(\theta)= 2^n a_n \cos^{2n}\frac{\theta}{2}\quad\text{and}\quad S(\theta) =-\sin\theta,
\end{equation}
where their $L^\infty$-norms satisfy \eqref{C-3}. It follows from \eqref{A-1} that 
\begin{equation} \label{C-3-0}
|\dot{\theta}_i-\nu_i|
\leq\frac{\kappa}{N} |S(\theta_i)| \sum_{j=1}^N|I_n(\theta_j)|\leq \kappa 2^n a_n.
\end{equation}
Hence, we have
\begin{equation} \label{C-3-1}
\dot{\theta}_i(t)\in\left(\nu_i- \kappa 2^n a_n,~\nu_i+ \kappa 2^n a_n \right), \quad i \in [N],~~t >0.
\end{equation}
Suppose that 
\[ \nu_i > \nu_j. \]
Then, it follows from \eqref{C-3-1} that 
\[ |{\dot \theta}_i(t) - {\dot \theta}_j(t)| \geq |\nu_i - \nu_j| - \kappa 2^{n+1} a_n  = \omega_{ij}^{\infty} > 0, \quad t \geq 0. \]
We take an infimum over $t$ to find the desired estimate. 
\end{proof}
\begin{remark} \label{R3.2}
Below, we briefly comment on the result of Theorem \ref{T3.1}. \newline
\begin{enumerate}
\item
Suppose that all natural frequencies are completely distributed:
\[ i \neq j \quad \Longrightarrow \quad \nu_i \neq \nu_j. \]
If the coupling strength satisfies 
\[ 0 \leq \kappa <  \kappa_{\mathrm{inc}},  \]
then we have
\[  \inf_{0 \leq t < \infty}  \min_{i \neq j} |\dot{\theta}_i(t)-\dot{\theta}_j(t)| \geq \min_{i \neq j} |\nu_i - \nu_j| - \kappa 2^{n+1} a_n > 0. \]
\item
By Lemma \ref{L3.1}, we have
\[ 2^{n +1}  a_n = {\color{black}{\mathcal O}(\sqrt{n})}, \quad n \gg 1. \]
Thus, for $n \gg 1$, the upper bound for $\kappa$ satisfies 
\[ \kappa^{ij}_{\mathrm{inc}} = {\mathcal O}\left(\frac{1}{\sqrt{n}}\right)|\nu_i - \nu_j|, \quad n \gg 1. \]
\item
 For identical ensemble with $\nu_i = \nu$, \eqref{C-3-0} implies
\[ |{\dot \theta}_i - \nu | \leq \kappa {\color{black}2^n}a_n \sim  \kappa\mathcal{O}(\sqrt{n}), \quad n \gg 1. \]
Again, this yields
{\color{black}\[  |\rho_i - \nu | \leq {\mathcal O}(\kappa \sqrt{n}), \]}
if $\rho_i$ exists. 
\end{enumerate}
In all cases, the coupling strength $\kappa_{\mathrm{inc}}$ does not depend on the number of oscillators. 
\end{remark}
\section{Emergence of death} \label{sec:4}
\setcounter{equation}{0}
In this section, we study sufficient frameworks leading to emergence of oscillator death state. Our strategy is to find a bounded invariant set of phases so that the corresponding rotation number becomes zero if they exist. 

\subsection{Complete oscillator death} \label{sec:4.1}
For $\alpha_n \in (0, \pi)$, we define a square box $ B(\alpha_n)$ in $\bbr^N$ as follows.
\begin{equation} \label{D-1}
B(\alpha_n) := \Big \{\Theta \in \bbr^N:~|\theta_i| < \alpha_n ~\text{ for all }~ i \in [N] \Big \}.
\end{equation}
In next lemma, we study the positive invariance of $B(\alpha_n)$. 
\begin{lemma}\label{L4.1} \emph{(Existence of a positive-invariant set)} 
Suppose parameters satisfy 
\begin{equation} \label{D-1-1}
0 < \alpha_n < \pi  \quad \mbox{and} \quad  \kappa> \frac{ \| {\mathcal V} \|_{\infty}}{|(SI_n)(\alpha_n)|} =: {\color{black}\kappa_{d,n}(\alpha_n)}.
\end{equation}
Then, $B(\alpha_n)$ in \eqref{D-1} is positive-invariant under \eqref{A-1} and \eqref{A-2}. 
\end{lemma}
\begin{proof}  Let $\Theta = \Theta(t)$ be a solution to \eqref{A-1} with initial data $\Theta^0$ satisfying
\begin{equation} \label{D-2}
\Theta^0 \in B(\alpha_n). 
\end{equation}
Then, it suffices to check that 
\begin{equation} \label{D-3}
 \Theta(t) \in  B(\alpha_n), \quad t > 0. 
\end{equation}
{\it Proof of \eqref{D-3}}: we use the continuity argument. For this, we define a set ${\mathcal T}$ and its supremum:
\[\mathcal{T}:=\{ t \in (0,\infty):\Theta(s)\in B(\alpha_n) ~\text{ for }~ s \in (0,t)\}\quad\mbox{and}\quad t^\infty:=\sup\mathcal{T}. \]
By \eqref{D-2} and the continuity of solution,  ${\mathcal T}$ contains a nontrivial interval. Hence, we have
\[  0 < t^{\infty} \leq \infty. \]

We claim $t^{\infty} = \infty$.
Suppose the contrary holds, $t^{\infty} < \infty$. Then, there exists an index $i_0 \in [N]$ which satisfies
\begin{equation} \label{D-4}
|\theta_{i_0}(t^\infty)|=\alpha_n, \quad \frac{d}{dt} |\theta_{i_0}| \Big|_{t=t^\infty} \geq 0.
\end{equation}
Moreover, we have
\[|\theta_j(t^\infty)|\le\alpha_n,\quad\forall j\ne i_0.\]
For such $i_0$, we use structural assumptions \eqref{B-1}, \eqref{B-2} and \eqref{D-1-1} to find 
\begin{align*}
\begin{aligned}
 \frac{d}{dt} |\theta_{i_0}| \Big|_{t=t^\infty} 
&=\nu_{i_0} ~\text{sgn}\left(\theta_{i_0}(t^{\infty})\right)+\frac{\kappa}{N}S(|\theta_{i_0}(t^{\infty})|)\sum_{j=1}^NI_n(\theta_j(t^{\infty}))\\
&\leq \| {\mathcal V} \|_{\infty} +\kappa (S I_n) (\alpha_n) <0.
\end{aligned}
\end{align*}
This is contradictory to $\eqref{D-4}_2$.  Therefore $t^{\infty} = \infty$ and we have the desired estimate. 
\end{proof}

\begin{remark} Lemma \ref{L4.1} already yields 
\[ \rho_i = 0, \quad i \in [N], \]
when $\Theta(t)$ is in $B(\alpha_n)$ once. This implies the emergence of complete oscillator death.  
\end{remark}

\begin{figure}
		\centering
			\includegraphics[width=0.6\textwidth]{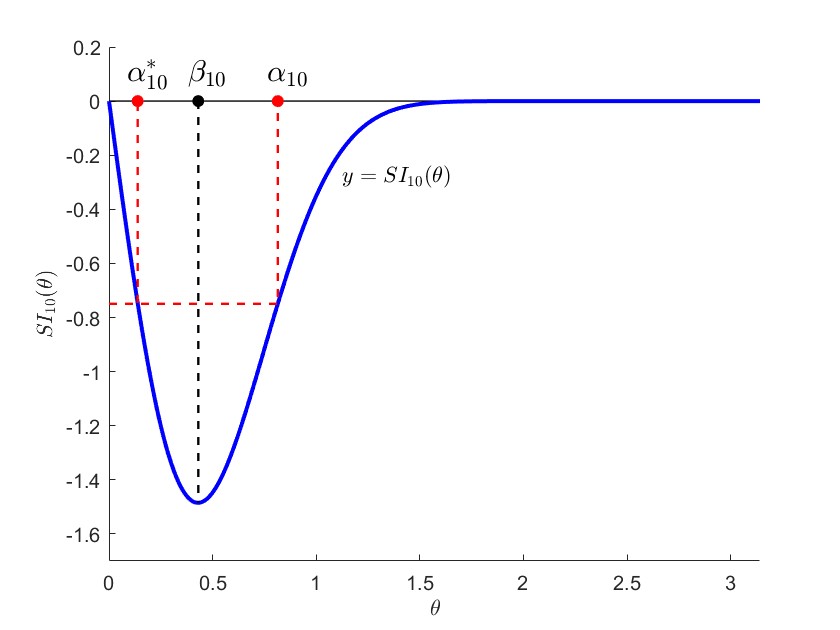} 
			\caption{Schematic diagrams of $SI_n$ with $n=10$} 	
		\label{Fig-2}
	\end{figure}

Figure \ref{Fig-2} shows the shape of $SI_n$, for the case of $n=10$, which determines $\kappa_{d,n}(\alpha_n)$ by \eqref{D-1-1}.
Let $\beta_n$ be the minimum point of the function $SI_n$ for $\theta \in (0,\pi)$: 
\begin{equation} \label{D-4-1}
\beta_n := \mbox{argmin}_{0 < \theta < \pi} (SI_n)(\theta). 
\end{equation}
By direct calculation, we have
\[ 
\beta_n = \cos^{-1} \Big( \frac{n}{n+ 1}  \Big ). \]
Then, it is easy to see that 
\begin{align*}
(SI_n)(\beta_n) \leq (SI_n)(\theta) \leq 0, \quad  \theta \in [0, \pi]. 
\end{align*}

\begin{remark} 
The minimal value $(SI_n)(\beta_n)$ is already computed in Lemma \ref{L3.1}, which is in the order of ${\mathcal O}(1)$ with respect to $n$. Hence, the value $\kappa_{d,n}(\alpha_n)$ can be set as in the same order ${\mathcal O}(1)$:
\[ \kappa_{d,n}(\alpha_n) \geq \kappa_{d,n} (\beta_n) \sim {\mathcal O}(1) \| {\mathcal V} \|_\infty, \quad n \gg 1. \]
\end{remark}

For $\alpha_n \in (\beta_n, \pi)$, let $\alpha_n^* \in (0, \pi)$ be the unique point satisfying 
\begin{align*}
\alpha_n^* \leq \beta_n \leq \alpha_n, \quad (SI_n)(\alpha_n^*) = (SI_n)(\alpha_n).
\end{align*}
In Lemma \ref{L4.1}, we have shown that the set  $B(\alpha_n)$  is positive-invariant. Next, we show that the set $B(\alpha_n^*) (\subset B(\alpha_n))$ further attracts all the trajectories issued from $B(\alpha_n)$ in finite time. 
\begin{theorem}\label{T4.1} 
\emph{(Complete oscillator death)}
Let $\beta_n \in (0, \pi)$ be the minimum point determined by \eqref{D-4-1}, and suppose parameters $\alpha_n, \kappa$ and initial data satisfy
\begin{align*}
\beta_n < \alpha_n < \pi,  \quad  \Theta^0 \in B(\alpha_n), \quad  \kappa>   \frac{\| {\mathcal V} \|_{\infty}}{|(SI_n)(\alpha_n)|} = \kappa_{d,n}(\alpha_n), 
\end{align*}
 and let $\Theta(t)$ be a global solution to \eqref{A-1}. Then the set $B(\alpha^*_n)$ attracts the trajectory $\{ \Theta(t) \}$ in finite-time, i.e., there exists a nonnegative constant $t_*$ such that
\[\Theta(t)\in B(\alpha^*_n), \quad t \geq t_*.\]
\end{theorem}
\begin{proof} By Lemma \ref{L4.1}, we have
\[ |\theta_j(t)| < \alpha_n \quad\text{for all}\quad t \geq 0. \]
Then, we use the structural condition of $S$ to see that for any $i$ with $|\theta_i(t)| \in (\alpha^*_n,\alpha_n)$,
\begin{align*}
\begin{aligned}
\frac{d|\theta_i(t)|}{dt}
&=\nu_i ~\text{sgn}\left(\theta_i(t)\right)+\frac{\kappa}{N}S(|\theta_i(t)|)\sum_{j=1}^NI_n(\theta_j(t))\\
&\le \| {\mathcal V} \|_{\infty} +\kappa S(|\theta_i(t)|)I_n(\alpha_n)\\
&\le\kappa\Big (-S(\alpha_n) +S(|\theta_i(t)|) \Big) I_n(\alpha_n)   <0.
\end{aligned}
\end{align*}
This implies that $|\theta_{i}(t)|$ should enter the interval $[0,\alpha^*_n)$ after some finite time $t_* \geq 0$.
\end{proof}
\begin{remark} Below, we comment on the results of Theorem \ref{T4.1}. 
\begin{enumerate}
\item
Note that the coupling condition $\kappa > \kappa_{d,n}(\alpha_n)$ is not only sufficient but also necessary to make both $B(\alpha_n)$ and $B(\alpha_n^*)$ positive-invariant. Existence of a positive-invariant set in the phase space could have technically complicated conditions which varies on the distribution of natural frequencies and initial data. However, for generic initial phase data with closely accumulated natural frequencies, one can expect that the condition for death state will mainly follow $\kappa >\kappa_{d,n}(\beta_n) \sim {\mathcal O}(1)\| {\mathcal V} \|_\infty$. This is also related to Theorem \ref{T4.2} and Remark \ref{R4.4} below.
\vspace{.1cm}
\item
One may also want to set a fixed bound $\alpha$ to make $B(\alpha)$ positive-invariant. In this case, $\alpha_n\equiv \alpha$ for all $n$. Since $\alpha>0$,
\[|(SI_n)(\alpha)|=a_n\sin\alpha(1+\cos\alpha)^n\sim\mathcal{O}\left(\sin^n\frac{\alpha}{2}\right)\]
which leads to
\[\kappa_{n,d}(\alpha)\sim\mathcal{O}(\gamma^n)\|\mathcal{V}\|_\infty\quad\mbox{for some}~~\gamma>1.\]
\end{enumerate}
\end{remark}

\subsection{Partial oscillator death} \label{sec:4.2}  In this subsection, we deal with sufficient conditions for partial oscillator death.  For a positive integer $2 \leq p \leq N$, we define a cylinder set as
\[ B_p(\alpha) : = \Big \{ \Theta = (\theta_1, \ldots, \theta_N) \in {\mathbb R}^N:~|\theta_i(t)|<\alpha_n,\quad \forall~ i \in [p]  \Big \}. \]
For $p  = N$, this cylinder coincides with $B(\alpha)$ defined in \eqref{D-1}. 
\begin{theorem}\label{T4.2}
Suppose parameters and initial data satisfy 
\begin{equation}  \label{D-6-1}
0 < \alpha_n < \pi, \quad  \max_{1 \leq i \leq p}  |\theta_i^0|<\alpha_n, \quad  \kappa>-\frac{N}{p}\frac{\| {\mathcal V} \|_{\infty}}{(SI_n)(\alpha_n)}, 
\end{equation}
and let $\Theta$ be a global solution to \eqref{A-1}--\eqref{A-2}. Then the following assertions hold.
\begin{enumerate}
\item
$B_p(\alpha)$ is positive-invariant. 
\vspace{0.2cm}
\item
For $i \in [N] \setminus [p]$, if there exists $t_* \geq 0$ such that $|\theta_i(t_*)|<\alpha_n,$ then we have 
\[  |\theta_i(t)|<\alpha_n\quad\mbox{for all }t\ge t_*.\]
\end{enumerate}
\end{theorem}
\begin{proof} Since the second assertion can be verified using the same argument as the first assertion, we first focus on the first statement below. For this, we again use continuity argument. By setting 
\[\mathcal{T}:= \Big \{\tau \in[0,\infty): |\theta_j(t)| \leq \alpha_n \quad \forall~ j \in [p], \quad \forall~t \in [0,\tau) \Big \}\quad\mbox{and}\quad \tau^\infty:=\sup\mathcal{T},\]
we have nonempty ${\mathcal T}$ and $0 < \tau^{\infty} \leq \infty$ from the initial condition and the continuity of solutions.
Now we claim:
\[ \tau^{\infty} = \infty. \]
Suppose the contrary holds, $\tau^{\infty} < \infty$. This implies that there exists an index $i_0$ with
\begin{equation} \label{D-7}
|\theta_{i_0}(\tau^\infty)| = \alpha_n \quad \mbox{and} \quad  \frac{d |\theta_{i_0}(t)|}{dt} \Big|_{t = \tau^{\infty}} \geq 0.
\end{equation}
For $t \in [0,\tau^\infty]$, one has 
\begin{align}
\begin{aligned} \label{D-8}
\frac{d|\theta_{i_0}(t)|}{dt}
&=\nu_{i_0} ~\text{sgn}\left(\theta_{i_0}(t)\right)+\frac{\kappa}{N}S(|\theta_{i_0}(t)|)\sum_{j=1}^NI_n(\theta_j(t))\\
&\leq\nu_{i_0} ~\text{sgn}\left(\theta_{i_0}(t)\right)+\frac{\kappa}{N}S(|\theta_{i_0}(t)|)\sum_{j=1}^p I_n(\theta_j(t))\\
&\leq  \| {\mathcal V} \|_{\infty} +\frac{\kappa p}{N} S(|\theta_{i_0}(t)|)I_n(\alpha_n).
\end{aligned}
\end{align}
Finally, we use $\eqref{D-6-1}_3$ and \eqref{D-8} to get 
\[
\frac{d|\theta_{i_0}(t)|}{dt}\Big|_{t=\tau^\infty} \leq \| {\mathcal V} \|_{\infty} +\frac{\kappa p}{N} S(\alpha_n)I_n(\alpha_n) < 0,
\]
which is contradictory to $\eqref{D-7}_2$. Hence, we have
\[  \tau^{\infty} = \infty \quad \mbox{and} \quad  |\theta_i(t)|<\alpha_n,\quad \forall~ i \in [p].  \]
\end{proof}

\begin{remark} \label{R4.4}
Note that Theorem \ref{T4.2} seems similar to Theorem \ref{T4.1}. However, it implicitly explains the emergence of convergence from more general initial data as follows.
\begin{enumerate}
\item
If oscillators are initially spreaded to $(-\pi,\pi)$ uniformly, then half of oscillators will be commonly in a half circle, $(-\pi/2,\pi/2)$. We set
\[  \alpha_n = \frac{\pi}{2}, \quad  (SI_n)(\alpha_n) = -a_n. \]
Therefore, Theorem \ref{T4.2} implies
\[ \kappa > \frac{2 \| {\mathcal V} \|_{\infty}}{a_n}\sim\mathcal{O}\left(\frac{\sqrt{n}}{2^n}\right)\|\mathcal{V}\|_\infty \quad \Longrightarrow \quad |\theta_i(t)|<\frac{\pi}{2},\quad \forall i \in [N]. \] 
Once there exists an invariant set, then $\Theta(t)$ is bounded.
\vspace{.1cm}
\item
As in Theorem \ref{T4.1}, once oscillators are trapped in the invariant set $[-\alpha_n,\alpha_n]$, they will be attracted to the set $[-\alpha^*_n,\alpha^*_n]$ in a finite time. 
\end{enumerate}
\end{remark}
\section{Emergence of locking}\label{sec:5}
\setcounter{equation}{0}
In this section, we study the emergence of phase locking for a Winfree ensemble with the same natural frequency. Here we mainly follow the methodology in \cite{H-K-P-R} and carefully extend it to the case with higher-order couplings. \newline

Let $\Theta = \Theta(t)$ be a global solution to \eqref{A-1}--\eqref{A-2}. Then, for $t \geq 0$, we set time-dependent extremal indices $M_t,~m_t$ and functionals $A(t), R(t)$:
\[
\begin{cases}
\displaystyle M_t:=\text{argmax}_i\theta_i(t),\quad m_t:=\text{argmin}_i\theta_i(t), \\
\displaystyle A(t):=\frac{\theta_{M_t}(t)+\theta_{m_t}(t)}{2},\quad R(t):=\frac{\theta_{M_t}(t)-\theta_{m_t}(t)}{2}.
\end{cases}
\]
Throughout this section, we assume that all oscillators have the same natural frequency:
\[  \nu_i = \nu, \quad \forall~i \in [N]. \]
This condition can be loosen as one can check results in \cite{H-K-P-R} with non-identical natural frequencies.

First, from the identical natural frequencies and the uniqueness of solutions, the order of oscillators does not change over time:
\[ \theta^0_i \le \theta_j^0,\quad \Longrightarrow \quad\theta_i(t)\le\theta_j(t),\quad \forall~t\ge0.\]
Therefore, maximal and minimal oscillators $\theta_M$ and $\theta_m$ are determined by that of initial phase oscillators. In the sequel, we consider the situation in which all oscillators are confined in a half circle geometrically:
\[ |\theta_M(t)-\theta_m(t)| < \pi, \]
so that 
\[ R(t) < \frac{\pi}{2}. \]
Of course, this a priori condition will be justified later. Moreover, we also consider a sufficient condition to make sure the strict increment of $A(t)$. Once the functional $A$ is strictly increasing, there exist an increasing sequence of times $\{ t_l^\pm \}$ such that
	\[A(t_l^-)=2l\pi-\frac{\pi}{2},\quad A(t_l^+)=2l\pi+\frac{\pi}{2}, \quad \l\in\mathbb{Z}. \]
	Without loss of generality, we set $N_*$ to be the smallest positive integer satisfying the following relation:
\begin{equation}\label{E-1}
0\leq t_{N_*}^-<t_{N_*}^+<t_{N_*+1}^-<t_{N_*+1}^+<\cdots.
\end{equation}
Then, the value $t_{N_*-1}^+$ can be positive or negative, but $t_{N_*-1}^- \leq 0$. 
\begin{figure}
		\centering
		\begin{subfigure}[ht]{0.45\textwidth}
\centering
			\includegraphics[width=0.95\textwidth]{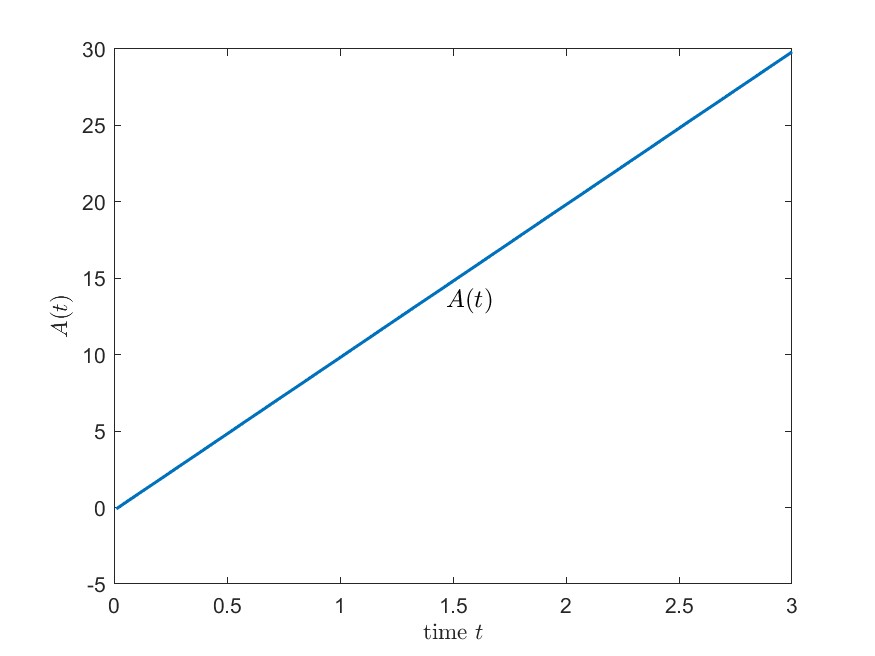}\caption{Temporal evolution of $A_t$}
		\end{subfigure}
\quad
		\begin{subfigure}[ht]{0.45\textwidth}
\centering
			\includegraphics[width=0.95\textwidth]{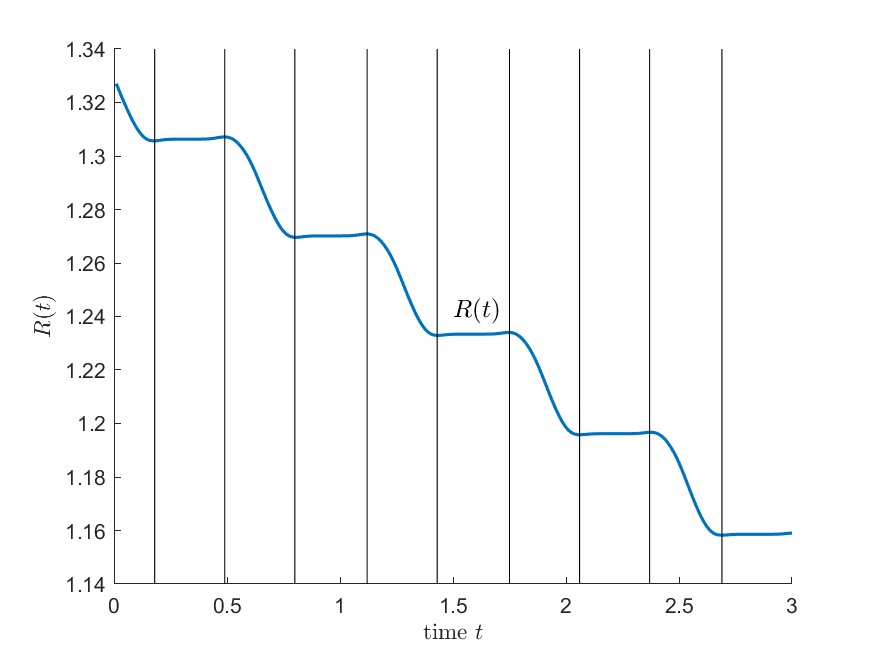}\caption{Temporal evolution of $R_t$ } 	
			\end{subfigure}
	\caption{Schematic diagrams of $A_t$ and $R_t$ under the framework of Theorem \ref{T5.1}}
		\label{Fig-3}
	\end{figure}

\subsection{Preparatory lemmas} \label{sec:5.1} 
In what follows, we present a series of basic estimates on $A$ and $R$. 
\begin{lemma}\label{L5.1}
Suppose that system parameters satisfy
\begin{equation} \label{E-2}
\nu_i = \nu > 0,~\forall~i \in [N] \quad \mbox{and} \quad  0 < \kappa  <\frac{\nu}{2^n a_n },
 \end{equation}
and let $\Theta=(\theta_1,\cdots,\theta_N)$ be a global solution to \eqref{A-1}--\eqref{A-2}. Then, functionals $A$ and $R$ satisfy the following three assertions:
\begin{enumerate}
\item
$A(t)$ is strictly increasing over time.
\item
If $t\in(t_l^-,t_l^+)$ and $0<R(t)\leq{\pi}/2$, then $R$ is strictly decreasing over time.
\item
If $t\in(t_l^+,t_{l+1}^-)$ and $0<R(t)\leq{\pi}/2$, then $R$ is strictly increasing over time.
\end{enumerate}		
\end{lemma}
\begin{proof}
(i)~We use the explicit form of $S$ in $\eqref{A-2}_2$:
		\[ \dot{\theta_i} = \nu +\frac{\kappa}{N}\sum_{j=1}^NS(\theta_i)I_n(\theta_j) = \nu -\frac{\kappa}{N}\sum_{j=1}^N\sin \theta_i I_n(\theta_j),\quad i \in [N]. \]
Then, with $\eqref{E-2}_2$, we get the desired increasing property of $A$:
\begin{align}
\begin{aligned} \label{E-2-0}
\dot{A} &=\frac{1}{2}(\dot \theta_{M} + \dot \theta_{m})
=\nu-\frac{\kappa}{2N}(\sin(\theta_M)+\sin(\theta_m))\sum_{j=1}^NI_n(\theta_j)  \\
&=  \nu - \kappa I_{n,c} \sin A \cos R \geq \nu - \kappa 2^n a_n > 0.
\end{aligned}
\end{align}
Here we used simplified notation $I_{n,c} := \sum_{j=1}^NI_n(\theta_j)$ and $\| I_n \|_{L^{\infty}} \leq 2^n a_n$ from \eqref{C-3-0-1}. \newline

(ii) and (iii):~From the definition of $R(t)$, one has
\begin{align}
\begin{aligned} \label{E-2-1}
\dot{R}(t) &=\frac{1}{2}(\dot \theta_{M}(t) - \dot \theta_{m}(t))
=\left(\frac{\kappa}{2N}\sum_{j=1}^NI_n(\theta_j(t))\right)(S(\theta_M(t))-S(\theta_m(t))) \\
&=-\kappa I_{n, c}(t) \cos A(t) \sin R(t).
\end{aligned}
\end{align}
If $t\in(t_l^-,t_l^+)$ and $0<R(t)\leq{\pi}/2$, we have
\[
\cos A(t) > 0\quad\mbox{and} \quad \sin R(t)>0.
\]
Then, \eqref{E-2-1} implies
\[ \dot{R}(t) < 0, \quad \mbox{i.e.,} \quad \mbox{$R$ is strictly decreasing over time}. \]
The remaining case, $t\in(t_l^+,t_{l+1}^-)$ and $0<R(t)\leq{\pi}/2$, can be treated similarly. 
\end{proof}
Note that Lemma \ref{L5.1} allows us to differentiate $R$ with respect to $A$.  For a given $\alpha >0$, we define functionals $L_1(A)$ and $L_2(A)$:
		\begin{align*}
		\begin{cases}
		\displaystyle L_1(A):=\frac{-\kappa\cos A}{\nu- \kappa 2^n a_n}\left(  \alpha  a_n \sqrt{2n-1} \left(\frac{2n-1}{n}\right)^{n-1} + I_n(A)\right),\\
		\displaystyle L_2(A):=\frac{\kappa\cos A}{\nu + \kappa 2^n a_n}\left(  \alpha  a_n \sqrt{2n-1} \left(\frac{2n-1}{n}\right)^{n-1} - I_n(A)\right).
	       \end{cases}
		\end{align*}
	
	\begin{lemma}\label{L5.2}
		Suppose that parameters satisfy
		\[ \nu_i = \nu > 0,~\forall~i \in [N], \quad  l\ge N_*, \quad \alpha\in(0,\pi/2], \quad  0 < \kappa <\frac{\nu}{2^n a_n},\]
		and  let $\Theta=(\theta_1,\cdots,\theta_N)$ be a global solution to \eqref{A-1}--\eqref{A-2}. Then, the following assertions hold.
		\begin{enumerate}
		\item
	         If $t\in(t_l^-,t_l^+)$ and $0<R(t)\le\alpha$, then we have
		\[L_1(A)\le\frac{1}{\sin R}\frac{dR}{dA}\le L_2(A).\]
		\item
		If $t\in(t_l^+,t_{l+1}^-)$ and $0<R(t)\le\alpha$, then we have
		\[L_2(A)\le\frac{1}{\sin R}\frac{dR}{dA}\le L_1(A). \]
		\end{enumerate}
		\end{lemma}
\begin{proof} It follows from \eqref{E-2-0} and \eqref{E-2-1} that 
\[ \dot A=\nu -\kappa I_{n,c} \sin A\cos R \quad\text{and}\quad \dot R=-\kappa I_{n,c} \cos A\sin R. \]
These yield
\begin{equation} \label{E-3-1}
\frac{1}{\sin R}\frac{dR}{dA}=\frac{-\kappa I_{n,c} \cos A  }{\nu -\kappa I_{n,c}\sin A\cos R}. 
\end{equation}
Note that  $\sin R$ is always positive since we assumed $R(t)\in (0,\pi/2]$. \newline

\noindent (i)~Suppose that $t\in(t_l^-,t_l^+)$ and  $0<R(t)\le\alpha$. Then, we have
		\[2l\pi-\frac{\pi}{2}<A(t)<2l\pi+\frac{\pi}{2}, \quad \mbox{i.e.,} \quad  \cos A>0.  \]
		It follows from $|I_{n,c}|\leq 2^n a_n$ that 
		\begin{equation}  \label{E-3-2}
		0<\nu- \kappa 2^n a_n < \nu - \kappa I_{n,c}\sin A\cos R < \nu  + \kappa  2^n a_n.
		\end{equation}
		Next, we estimate the average influence $I_{n,c}$ with $I_n(A)$. First, note that
		\[ |I_n(\theta_j) - I_n(A)| \leq \|I_n'\|_{\infty} |\theta_j-A|  \le a_n \sqrt{2n-1} \left(\frac{2n-1}{n}\right)^{n-1}|\theta_j-A|. \]			
		Under the assumption $0<R(t) \le \alpha$, the above estimate and 
		\[ 
		|\theta_j-A|  \leq \frac{1}{2} \Big(  | \theta_j - \theta_M | + |\theta_j - \theta_m| \Big) \leq R(t)
		\]
		induces the estimation:
		\begin{align}
		\begin{aligned} \label{E-3-3}
		|I_{n,c}(\Theta) - I_n(A)|  &\leq \frac{1}{N} \sum_{j=1}^{N} |I_n(\theta_j) - I_n(A) | \leq a_n \sqrt{2n-1} \left(\frac{2n-1}{n}\right)^{n-1} R(t) \\
		&\le \alpha  a_n \sqrt{2n-1} \left(\frac{2n-1}{n}\right)^{n-1}.
		\end{aligned}
		\end{align}
		In \eqref{E-3-1}, we combine \eqref{E-3-2} and \eqref{E-3-3} to find 
\[L_1(A)\le\frac{1}{\sin R}\frac{dR}{dA}\le L_2(A).\] 

\noindent (ii)~Suppose that $t\in(t_l^-,t_{l+1}^+)$ and $0<R(t)\le\alpha$. Then, we have
		\[2l\pi+\frac{\pi}{2}<A(t)<2(l+1)\pi-\frac{\pi}{2},\quad\mbox{i.e.,}\quad\cos A<0.\]
		Therefore, the opposite inequalities hold due to the sign of  $\cos A$.
	\end{proof}
Next, we show that  the functional $R(t)$ is bounded. In fact,  for identical ensemble, we may eventually prove that $R(t)$ goes to zero as $t\to\infty$. Since the case for $n=1$ was already studied in \cite{H-K-P-R}, we assume $n\ge2$. Lemma \ref{L5.1} shows that $R(t)$ increases for $t\in(t_l^+,t_{l+1}^-)$ and decreases for $t\in(t_l^-,t_l^+)$. From Lemma \ref{L5.2}, we need to estimate, first, how small $R(t)$ increases by measuring the integral of $L_2(A)$, second, how deeply $R(t)$ decreases from the integral of $L_1(A)$. This comparison needs the following computation. 
\begin{lemma}\label{L5.3}
For $n\ge2$, the following estimates hold.
\begin{align}
\begin{aligned} \label{E-4}
& (i)~\int_{\frac{\pi}{2}}^{\frac{3\pi}{2}}(-\cos A)\cos^{2n}\frac{A}{2}dA \leq  \frac{1}{2^{n-1}}. \\
& (ii)~\int_{-\frac{\pi}{2}}^{\frac{\pi}{2}}\cos A \cos^{2n}\frac{A}{2}dA \geq \frac{n}{n+1}\frac{\pi+2}{{\color{black}2^n}a_n} \sim {\mathcal O}\left(\frac{1}{\sqrt{n}}\right). \\
& (iii)~\int_{\frac{\pi}{2}}^{\frac{3\pi}{2}}(-\cos A)\cos^{2n}\frac{A}{2}dA {\color{black}  = \int_{-\frac{\pi}{2}}^{\frac{\pi}{2}}\cos A \cos^{2n}\frac{A}{2}dA - \frac{\pi}{2^na_n}\frac{2n}{n+1}}.
\end{aligned}
\end{align}
\end{lemma}

\begin{proof} 
(i)~The first estimation is from the maximal value of $\cos^{2n}\frac{A}{2}$:
\begin{align*}
\int_{\frac{\pi}{2}}^{\frac{3\pi}{2}}(-\cos A)\cos^{2n}\frac{A}{2}dA &\leq \int_{\frac{\pi}{2}}^{\frac{3\pi}{2}}(-\cos A)\frac{1}{2^n}dA = \frac{1}{2^{n-1}}.
\end{align*}

\vspace{0.2cm}

\noindent (ii)~For any $n\ge2$, we use a well-known trigonometric formula
		\[ \int\cos^n \theta d\theta =\frac{1}{n}\cos^{n-1} \theta \sin \theta +\frac{n-1}{n}\int\cos^{n-2} \theta d\theta \]
recursively to see
\begin{align}
\begin{aligned} \label{E-5}
& \int_{-\frac{\pi}{4}}^{\frac{\pi}{4}}\cos^{2n}\theta d\theta = \frac{1}{n} \left( \cos^{2n-1}\frac{\pi}{4}\sin\frac{\pi}{4} \right) + \frac{2n-1}{2n} \int_{-\frac{\pi}{4}}^{\frac{\pi}{4}} \cos^{2n-2}\theta d\theta \\
&\hspace{0.2cm}=\frac{1}{n}\left(\frac{1}{2}\right)^n +\sum_{k=1}^{n-1} \Big[ \frac{(2n-1)(2n-3)\cdots( 2n - (2k-1)))}{(2n)(2n-2)\cdots(2n - (2k-2)))}\cdot\frac{1}{n-k}\left(\frac{1}{2}\right)^{n-k} \Big ] \\
&\hspace{0.2cm} +\frac{(2n-1)(2n-3)\cdots (2n - (2n -1))}{(2n)(2n-2)\cdots (2n - (2n-2))}\int_{-\frac{\pi}{4}}^{\frac{\pi}{4}}\cos^{0}\theta d\theta,\\
&\hspace{0.2cm} =\frac{1}{n}\left(\frac{1}{2}\right)^n  +\sum_{k=1}^{n-1} \Big[ \frac{(2n-1)(2n-3)\cdots( 2n - (2k-1)))}{(2n)(2n-2)\cdots(2n - (2k-2)))}\cdot\frac{1}{n-k}\left(\frac{1}{2}\right)^{n-k} \Big ]  +\frac{\pi}{2^{n+1}a_n}.
\end{aligned}
\end{align}
Note that 
\begin{align*}
\begin{aligned}
&\int_{-\frac{\pi}{2}}^{\frac{\pi}{2}}\cos A \cos^{2n}\frac{A}{2}dA =\int_{-\frac{\pi}{2}}^{\frac{\pi}{2}}\left(2\cos^2\frac{\theta}{2}-1\right)\left(\cos^{2n}\frac{\theta}{2}\right)d\theta \\
& \hspace{1cm} =2\int_{-\frac{\pi}{4}}^{\frac{\pi}{4}}(2\cos^2\theta-1)(\cos^{2n}\theta)d\theta =4\int_{-\frac{\pi}{4}}^{\frac{\pi}{4}}\cos^{2n+2}\theta-2 \int_{-\frac{\pi}{4}}^{\frac{\pi}{4}}\cos^{2n}\theta d\theta.
\end{aligned}
\end{align*}
On the other hand, we use \eqref{E-5} to find
\begin{align*}
\begin{aligned}
&\int_{-\frac{\pi}{2}}^{\frac{\pi}{2}}\cos A \cos^{2n}\frac{A}{2}dA = 4\int_{-\frac{\pi}{4}}^{\frac{\pi}{4}}\cos^{2n+2}\theta d\theta-2\int_{-\frac{\pi}{4}}^{\frac{\pi}{4}}\cos^{2n}\theta d\theta\\
& \hspace{0.2cm} = 4\left[ \frac{1}{n+1}\left(\frac{1}{2}\right)^{n+1}+\sum_{k=1}^{n}\frac{(2n+1)(2n-1)\cdots(2n-2k + 3)}{(2n+2)(2n)\cdots(2n - 2k +4)}\cdot\frac{1}{ n+1-k}\left(\frac{1}{2}\right)^{n+1 -k}+\frac{\pi}{2^{n+2} a_{n+1}} \right] \\
&  \hspace{0.2cm}\quad -2\left[ \frac{1}{n}\left(\frac{1}{2}\right)^n+\sum_{k=1}^{n-1}\frac{(2n-1)(2n-3)\cdots(2n-(2k-1))}{(2n)(2n-2)\cdots(2n-(2k-2))}\cdot\frac{1}{n-k}\left(\frac{1}{2}\right)^{n-k} +\frac{\pi}{2^{n +1}a_n} \right]\\
&  \hspace{0.2cm} = \left(\frac{1}{2}\right)^n \left[ \frac{2}{n+1} - \frac{2}{n} + \frac{2n+1}{2n+2}\cdot\frac{4}{n} \right] + \frac{\pi}{2^n a_n} \left[\frac{2n+1}{n+1} - 1 \right] \\
&  \hspace{0.2cm}\quad+ \sum_{k=1}^{n-1}\left[ 4\cdot\frac{2n+1}{2n+2} - 2 \right] \frac{(2n-1)(2n-3)\cdots(2(k+1)-1)}{(2n)(2n-2)\cdots(2(k+1))}\cdot\frac{1}{k}\left(\frac{1}{2}\right)^k \\
& \hspace{0.2cm}= \left(\frac{1}{2}\right)^n \left[\frac{4}{n+1}\right] + \frac{\pi}{2^n a_n} \left[ \frac{n}{n+1} \right] 
+ \sum_{k=1}^{n-1}\left[ \frac{2n}{n+1} \right] \frac{(2n-1)(2n-3)\cdots(2(k+1)-1)}{(2n)(2n-2)\cdots(2(k+1))}\cdot\frac{1}{k}\left(\frac{1}{2}\right)^k. \\
\end{aligned}
\end{align*}
where we used the relation:
\[ \frac{\pi}{2^{n} a_{n+1} }= \frac{2\pi}{2^{n} a_n} \Big( \frac{2n + 1}{2n +2} \Big).  \]
By selecting the second term and the third term with $k=1$, we get a lower bound of the integration:
\begin{align*}
\int_{-\frac{\pi}{2}}^{\frac{\pi}{2}}\cos A \cos^{2n}\frac{A}{2}dA \geq \frac{n}{n+1}\frac{\pi}{2^n a_n} + \frac{2n}{n+1}\frac{2}{2^n a_n}\frac{1}{2} =\frac{n}{n+1}\frac{\pi+2}{2^n a_n}.
\end{align*}
%{\color{red} Dongnam Ko: Why this is here?\\
%For $n\ge2$,
%\begin{align*}
%\int_{-\frac{\pi}{2}}^{\frac{\pi}{2}}\cos A\cos^{2n}\frac{A}{2}dA\le\frac{1}{3}+\frac{\pi}{2^na_n}+\sum_{k=1}^{n-1}\frac{2}{k}\left(\frac{1}{2}\right)^k
%\end{align*}
%}

\vspace{0.2cm}

\noindent (iii)~
For an integer $n\ge2$, we use
		\[\int\sin^n\theta d\theta=-\frac{1}{n}\sin^{n-1}\theta\cos\theta+\frac{n-1}{n}\int\sin^{n-2}\theta d\theta \]
to see
\[
			\int_{-\frac{\pi}{4}}^{\frac{\pi}{4}}\sin^{2n}\theta d\theta=-\frac{1}{n}\left(\frac{1}{2}\right)^n-\sum_{k=1}^{n-1}\frac{(2n-1)\cdots(2(k+1)-1)}{(2n)\cdots(2(k+1))}\cdot\frac{1}{k}\left(\frac{1}{2}\right)^k+\frac{(2n-1)\cdots1}{(2n)\cdots2}\cdot\frac{\pi}{2}. \]
Hence, we have
{\color{black}
\begin{align*}
\begin{aligned}
\int_{\frac{\pi}{2}}^{\frac{3\pi}{2}}-\cos &A\cos^{2n}\frac{A}{2}dA = 2\int_{-\frac{\pi}{4}}^{\frac{\pi}{4}}\sin^{2n}\theta d\theta-4\int_{-\frac{\pi}{4}}^{\frac{\pi}{4}}\sin^{2n+2}\theta d\theta\\
&=\left(\frac{1}{2}\right)^n\frac{4}{n+1}-\frac{\pi}{2^na_n}\frac{n}{n+1}+\sum_{k=1}^{n-1}\frac{2n}{n+1}\frac{(2n-1)\cdots( 2(k+1)-1 )}{(2n)\cdots ( 2(k+1) )}\cdot\frac{1}{k}\left(\frac{1}{2}\right)^k,
\end{aligned}
\end{align*}
which deduces the desired estimate by comparing it with the estimate of $\int_{-\frac{\pi}{2}}^{\frac{\pi}{2}}\cos A \cos^{2n}\frac{A}{2}dA$.
}
\begin{comment}
\begin{align}
\begin{aligned} \label{E-7}
&\int_{-\frac{\pi}{2}}^{\frac{\pi}{2}}\cos A \cos^{2n}\frac{A}{2}dA = 4\int_{-\frac{\pi}{4}}^{\frac{\pi}{4}}\sin^{2n+2}\theta d\theta-2\int_{-\frac{\pi}{4}}^{\frac{\pi}{4}}\sin^{2n}\theta d\theta\\
& \hspace{1cm} = \left(\frac{1}{2}\right)^n \left[-\frac{4}{n+1}\right] + \frac{\pi}{a_n} \left[ \frac{n}{n+1} \right]  \\
& \hspace{1cm} + \sum_{k=1}^{n-1}\left[ \frac{-2n}{n+1} \right] \frac{(2n-1)(2n-3)\cdots(2(k+1)-1)}{(2n)(2n-2)\cdots(2(k+1))}\cdot\frac{1}{k}\left(\frac{1}{2}\right)^k \\
& \hspace{1cm} \leq \frac{n}{n+1}\frac{\pi}{2^n a_n} - \frac{2n}{n+1}\frac{2}{2^n a_n}\frac{1}{2} =\frac{n}{n+1}\frac{\pi-2}{2^n a_n}.
\end{aligned}
\end{align}

Finally, we combine \eqref{E-6} and \eqref{E-7} to find the desired estimate. 
\end{comment}
\end{proof}

%%%%%%%%%%%%%%%%%%%%%%%%%%%%%%%%%%%%%%%%%%%%%%%%%%%%%%%%%%%%%%%%%%%%%%%%%%%%%%%%%%%%%%%%%%%
%%%%%%%%%%%%%%%%%%%%%%%%%%%%%%%%%%%%%%%%%%%%%%%%%%%%%%%%%%%%%%%%%%%%%%%%%%%%%%%%%%%%%%%%%%%
%%%%%%%%%%%%%%%%%%%%%%%%%%%%%%%%%%%%%%%%%%%%%%%%%%%%%%%%%%%%%%%%%%%%%%%%%%%%%%%%%%%%%%%%%%%\\
%%%%%%%%%%%%%%%%%%%%%%%%%%%%%%%%%%%%%%%%%%%%%%%%%%%%%%%%%%%%%%%%%%%%%%%%%%%%%%%%%%%%%%%%%%%\\
%%%%%%%%%%%%%%%%%%%%%%%%%%%%%%%%%%%%%%%%%%%%%%%%%%%%%%%%%%%%%%%%%%%%%%%%%%%%%%%%%%%%%%%%%%%
Next, we combine Lemma \ref{L5.2} and Lemma \ref{L5.3} to derive sufficient conditions on $\kappa$ and $\alpha$ for the boundedness of $R(t)$.
\begin{lemma}\label{L5.4} 
Suppose that system parameters satisfy 
\begin{align*}
&n \geq 2, \qquad 0 < \kappa \le \frac{\nu}{2^{n+1} a_n} \sim {\mathcal O} \left(\frac{\nu}{\sqrt{n}}\right) \quad\text{and}\\
&\alpha<\left(\frac{\pi}{2^{n+1}a_n}\frac{n}{n+1}-\frac{1}{2^n}\right)\frac{1}{\sqrt{2n-1}}\left(\frac{2n}{2n-1}\right)^{n-1}\sim\mathcal{O}\left(\frac{1}{n}\right).
\end{align*}
Let $\Theta$ be a global solution to \eqref{A-1}--\eqref{A-2} satisfying a priori condition: 
\[ \sup_{0 \leq t < \infty} R(t) \leq \alpha < \frac{\pi}{2}. \]
Then, the following estimates hold. 
\begin{align*}
\begin{aligned}
& (i)~\int_{\frac{\pi}{2}}^{\frac{3\pi}{2}}L_1(A)dA \leq {\color{black}\alpha\sqrt{2n-1}\left(\frac{2n-1}{2n}\right)^{n-1}} + \frac{1}{2^{n-1}}. \\
& (ii)~\int_{-\frac{\pi}{2}}^{\frac{3\pi}{2}} \frac{1}{\sin R}\frac{dR}{dA} dA \leq  \int_{-\frac{\pi}{2}}^{\frac{\pi}{2}}L_2(A)dA + \int_{\frac{\pi}{2}}^{\frac{3\pi}{2}}L_1(A)dA < 0.
\end{aligned}
\end{align*}

\end{lemma}
\begin{proof}
(i)~We use the integration of $\cos A$ and Lemma \ref{L5.3} to find 
\begin{align}
\begin{aligned} \label{E-8}
\int_{\frac{\pi}{2}}^{\frac{3\pi}{2}}L_1(A)dA &= \int_{\frac{\pi}{2}}^{\frac{3\pi}{2}}\frac{-\kappa\cos A}{\nu - \kappa 2^n a_n}\Big[ \alpha 2^n a_n\frac{\sqrt{2n-1}}{2}\left(\frac{2n-1}{2n}\right)^{n-1} + I_n(A) \Big]dA \\
&= \frac{\kappa 2^n a_n}{\nu-\kappa 2^n a_n}\Big [ {\color{black}\alpha\sqrt{2n-1}}\left(\frac{2n-1}{2n}\right)^{n-1} + \int_{\frac{\pi}{2}}^{\frac{3\pi}{2}}(-\cos A)\cos^{2n}\frac{A}{2}dA \Big ] \\
&\leq \frac{\kappa 2^n a_n}{\nu- \kappa 2^n a_n}\Big [ {\color{black}\alpha\sqrt{2n-1}}\left(\frac{2n-1}{2n}\right)^{n-1} + \frac{1}{2^{n-1}} \Big],
\end{aligned}
\end{align}
where the last inequality is due to \eqref{E-4} (i) in Lemma \ref{L5.3}. On the other hand, since the assumption on $\kappa$ guarantees
\[ \frac{\kappa 2^n a_n}{\nu- \kappa 2^n a_n} \leq 1, \]
we have
\begin{align*}
& \int_{\frac{\pi}{2}}^{\frac{3\pi}{2}}L_1(A)dA \leq {\color{black}\alpha\sqrt{2n-1}}\left(\frac{2n-1}{2n}\right)^{n-1}  + \frac{1}{2^{n-1}}.
\end{align*}
(ii) We use Lemma \ref{L5.2} to get 
\begin{equation} \label{E-9}
\int_{-\frac{\pi}{2}}^{\frac{3\pi}{2}} \frac{1}{\sin R}\frac{dR}{dA} dA \leq  \int_{-\frac{\pi}{2}}^{\frac{\pi}{2}}L_2(A)dA + \int_{\frac{\pi}{2}}^{\frac{3\pi}{2}}L_1(A)dA.
\end{equation}
Then, the right-hand side of  \eqref{E-9} becomes
\begin{equation}\label{E-9-1}
\begin{aligned}
&\int_{-\frac{\pi}{2}}^{\frac{\pi}{2}}L_2(A)dA+\int_{\frac{\pi}{2}}^{\frac{3\pi}{2}}L_1(A)dA \\
&\quad= \int_{-\frac{\pi}{2}}^{\frac{\pi}{2}}\frac{\kappa\cos A}{\nu +\kappa 2^n a_n }\left(\alpha 2^n a_n\frac{\sqrt{2n-1}}{2}\left(\frac{2n-1}{2n}\right)^{n-1} - I_n(A)\right)dA + \int_{\frac{\pi}{2}}^{\frac{3\pi}{2}}L_1(A)dA  \\
&\quad= \frac{\kappa 2^n a_n}{\nu + \kappa 2^n a_n}\left({\color{black}\alpha\sqrt{2n-1}}\left(\frac{2n-1}{2n}\right)^{n-1} - \int_{-\frac{\pi}{2}}^{\frac{\pi}{2}}\cos A\cos^{2n}\frac{A}{2}dA \right) \\
&\quad\quad + \frac{\kappa 2^n a_n}{\nu - \kappa 2^n a_n}\left({\color{black}\alpha\sqrt{2n-1}}\left(\frac{2n-1}{2n}\right)^{n-1}+\int_{\frac{\pi}{2}}^{\frac{3\pi}{2}}(-\cos A)\cos^{2n}\frac{A}{2}dA \right).
\end{aligned}
\end{equation}
In order to compare two different coefficients, we use the assumption on $\kappa$:
\[  \kappa  \leq \frac{\nu}{2^{n+1} a_n}~\left(\mbox{hence, }~\frac{\kappa 2^na_n}{\nu-\kappa2^na_n}\le\frac{3\kappa2^na_n}{\nu+\kappa2^na_n}\right).\]
Then, the integral \eqref{E-9-1} can be estimated as
\begin{equation*}
\begin{aligned}
&\int_{-\frac{\pi}{2}}^{\frac{\pi}{2}}L_2(A)dA+\int_{\frac{\pi}{2}}^{\frac{3\pi}{2}}L_1(A)dA \\
&\le \frac{\kappa 2^n a_n}{\nu + \kappa 2^n a_n} \left[ {\color{black}4}\alpha \sqrt{2n-1}\left(\frac{2n-1}{2n}\right)^{n-1} - \int_{-\frac{\pi}{2}}^{\frac{\pi}{2}}\cos A\cos^{2n}\frac{A}{2}dA \right. 
 \left. + 3\int_{\frac{\pi}{2}}^{\frac{3\pi}{2}}(-\cos A)\cos^{2n}\frac{A}{2}dA \right].
\end{aligned}
\end{equation*}

Next, we erase the second term using \eqref{E-4}(iii) of Lemma \ref{L5.3} and estimate the last term using \eqref{E-4}(i) to get 
\begin{equation*}
\begin{aligned}
&\int_{-\frac{\pi}{2}}^{\frac{\pi}{2}}L_2(A)dA+\int_{\frac{\pi}{2}}^{\frac{3\pi}{2}}L_1(A)dA \\
&\quad =\frac{\kappa 2^n a_n}{\nu +\kappa 2^n a_n}\left[4\alpha\sqrt{2n-1}\left(\frac{2n-1}{2n}\right)^{n-1}-\frac{\pi}{2^na_n}\frac{2n}{n+1}+2\int_{\frac{\pi}{2}}^{\frac{3\pi}{2}}(-\cos A)\cos^{2n}\frac{A}{2}dA\right]\\
&\quad <\frac{\kappa 2^n a_n}{\nu +\kappa 2^n a_n}\left[4\alpha \sqrt{2n-1}\left(\frac{2n-1}{2n}\right)^{n-1}-\frac{2\pi}{2^na_n}\frac{n}{n+1}+\frac{1}{2^{n-2}}\right].
\end{aligned}
\end{equation*}
%Dongnam Ko: previous estimation does not satisfied when $n=2$.

Therefore, it becomes a negative value if 
\[\alpha<\left(\frac{\pi}{2^{n+1}a_n}\frac{n}{n+1}-\frac{1}{2^n}\right)\frac{1}{\sqrt{2n-1}}\left(\frac{2n}{2n-1}\right)^{n-1}.\]
From the bounds of $a_n$ in \eqref{C-2-1} of Remark \ref{R3.1},
$a_n \leq \sqrt{n}/2^{n-1}$,
the right-hand side is always positive for $n \geq 2$.
\end{proof}

Next, we define constants $C_1^+, C_1^-, C_2^+$, and $C_2^-$ as follows:
\begin{align*}
\begin{aligned}
&-C_1^+ := \int_{2l\pi-\frac{\pi}{2}}^{2l\pi+\frac{\pi}{2}}L_1(A)dA, \quad 
C_1^- := \int_{2l\pi+\frac{\pi}{2}}^{2l\pi+\frac{3\pi}{2}}L_1(A)dA,\\
&-C_2^+ := \int_{2l\pi-\frac{\pi}{2}}^{2l\pi+\frac{\pi}{2}}L_2(A)dA,\quad
C_2^- := \int_{2l\pi+\frac{\pi}{2}}^{2l\pi+\frac{3\pi}{2}}L_2(A)dA.
\end{aligned}		
\end{align*} 
%These constants are all positive from the following lemma.
\begin{lemma}\label{L5.5}
Suppose that parameters satisfy 
\begin{align*}
\begin{aligned}
& n \geq 2, \quad l \in \mathbb Z, \quad 0 <  \kappa < \frac{\nu}{2^{n+1} a_n} \sim {\mathcal O}\left(\frac{\nu}{\sqrt{n}}\right), \\
&\alpha<\left(\frac{\pi}{2^{n+1}a_n}\frac{n}{n+1}-\frac{1}{2^n}\right)\frac{1}{\sqrt{2n-1}}\left(\frac{2n}{2n-1}\right)^{n-1}\sim\mathcal{O}\left(\frac{1}{n}\right),
\end{aligned}
\end{align*}
and let $\Theta$ be a global solution to \eqref{A-1}--\eqref{A-2} satisfying a priori condition:
\[ \sup_{0 \leq t < \infty} R(t) \leq \alpha < \frac{\pi}{2}. \]
Then, $C_1^+, C_1^-, C_2^+$, and $C_2^-$ are all positive and the following estimates hold:
\begin{align*}
&(i)-C_1^+<\ln\frac{R(t_l^+)}{R(t_l^-)}<-\frac{\sin\alpha}{\alpha}C_2^+\quad\text{and}\quad \frac{\sin\alpha}{\alpha}C_2^-<\ln\frac{R(t_{l+1}^-)}{R(t_l^+)}<C_1^-,\\
&(ii)\ln\frac{R(t_{l+1}^-)}{R(t_l^-)}<\frac{\sin\alpha}{\alpha}(C_1^--C_2^+) < 0\quad\text{and}\quad\ln\frac{R(t_{l+1}^+)}{R(t_l^+)}\ge-C_1^++C_2^-.
\end{align*}
	\end{lemma}
	\begin{proof}
		The positivity of constants follows similar arguments to the proof of Lemma \ref{L5.4} using \eqref{E-8} and \eqref{E-9-1} with smallness condition on $\alpha$. The proof of remaining estimates needs direct but lengthy calculation from Lemmas \ref{L5.2} to \ref{L5.4}. Since this uses structurally the same argument as in Lemma 3.5 of \cite{H-K-P-R}, here we omit it.
	\end{proof}
From Lemma \ref{L5.5}, we have the following corollary, of which proof is the same as in Corollary 3.1 of \cite{H-K-P-R}.
\begin{corollary}
Assume the positive parameters $l,N_*,\alpha$ satisfy
\[l\ge N_*,\quad\alpha\in(0,\alpha^\infty),\quad\sup_{0\le t<\infty}R(t)\le\alpha.\]
Then, we have the following estimates:
\begin{align*}
	&(i)~~\ln\frac{R(t_l^-)}{R(t_{N_*}^-)}\le-\frac{\sin\alpha}{\alpha}\cdot\frac{(\Omega-2K)}{2\pi}|C_2^+-C_1^-|\cdot|t_l^--t_{N_*}^-|.\\
	&(ii)~~\ln\frac{R(t_l^+)}{R(t_{N_*}^+)}\ge-\frac{(\Omega+2K)}{2\pi}|C_1^+-C_2^-|\cdot|t_l^+-t_{N_*}^+|.
\end{align*}
\end{corollary}

\subsection{Exponential decay of phase diameter} \label{sec:5.2}	
In this subsection, we present our last main result on the exponential decay of $R$. 
\begin{theorem}\label{T5.1}
Suppose that system parameters and initial data satisfy
\begin{align*}
&0 < \kappa < \frac{\nu}{2^{n+1} a_n} \sim {\mathcal O}\left(\frac{\nu}{\sqrt{n}}\right), \quad 
				R(0)<\alpha\exp \Big[ -{\color{black}\alpha\sqrt{2n-1}}\left(\frac{2n-1}{2n}\right)^{n-1}-\frac{1}{2^{n-1}} \Big],\\
&0<\alpha<\left(\frac{\pi}{2^{n+1}a_n}\frac{n}{n+1}-\frac{1}{2^n}\right)\frac{1}{\sqrt{2n-1}}\left(\frac{2n}{2n-1}\right)^{n-1}\sim\mathcal{O}\left(\frac{1}{n}\right),
\end{align*}
and let $\Theta(t)$ be a global solution to \eqref{A-1}. Then, there exists positive constants ${\tilde \beta}_i$ and ${\tilde \Lambda}_i$, $i=1,2$ such that
		{\color{black}\[ {\tilde \beta}_1e^{- {\tilde \Lambda}_1t}R(0)\le R(t)\le {\tilde \beta}_2e^{-{\tilde \Lambda}_2t}R(0) < \alpha,\quad t\ge0.\]}
	\end{theorem}
\begin{proof} Firstly, we will show 
\[R(t)<\alpha,\quad\forall t\ge0.\]
Before we start, let $N_*$ be the positive integer defined by \eqref{E-1}, where $t_{N_*}^-> 0$ and $t_{N_*-1}^- \leq 0$ by definition. Define a set
\[\mathcal{T}:=\{T\ge0:R(t)<\alpha,\quad\forall t\in[0,T]\}\]
and let $T^*:=\sup\mathcal{T}$. Since $0\in\mathcal{T}$, $T^*$ is a positive number. What we have to show is equivalent to $T^*=\infty$. Let us use the contradiction argument. Suppose $T^*<\infty$, which means $R(T^*)=\alpha$. We begin with showing $T^*\ge t_{N_*}^-$. We consider two cases with the assumption
\[T^*\in(0,t_{N_*}^-).\]

\noindent $\bullet$~Case A ($t_{N_*-1}^+\le0$): Since $A(0)\le 0$ and $L_1(A)\ge0$, we have
\[\frac{1}{\sin R(t)}\frac{dR}{dA}\le L_1(A(t)),\quad\forall t\in[0,t_{N_*}^-].\]
So, one has
\begin{align}\label{5.15}
\begin{aligned}
	\int_{A(0)}^{A(T^*)}\frac{1}{\sin R(t)}\frac{dR}{dA}dA&\le\int_{A(0)}^{A(T^*)}\L_1(A(t))dA\\
	&\le\int_{A(t_{N_*-1}^+)}^{A(t_{N_*}^-)}L_1(A(t))dA=\int_{-\frac{3\pi}{2}}^{-\frac{\pi}{2}}L_1(A(t))dA=C_1^-,
\end{aligned}
\end{align}
where we used the fact $t_{N_*-1}^+\le0\le T^*\le T_{N_*}^-$.
On the other hand, consider a function $f$ defined by
\[f(x):=-\ln\left(\frac{1+\cos x}{\sin x}\right),\quad x\in\left(0,\frac{\pi}{2}\right),\]
hence,
\[f'(x)=\frac{1}{\sin x}>0.\]
Then, by simple calculations, we get
\begin{align}\label{formula_f}
	\frac{\sin y}{y}(f(y)-f(x))<\ln\frac{y}{x}<f(y)-f(x).
\end{align}
Back to our story, from \eqref{5.15} with \eqref{formula_f}, we can derive
\begin{align*}
	\ln\left(\frac{R(T^*)}{R(0)}\right)<f(R(T^*))-f(R(0))\le C_1^-,
\end{align*}
which implies $\frac{R(T^*)}{R(0)}\le e^{C_1^-}$, that is,
\[R(T^*)\le R(0)e^{C_1^-}\le R(0)e^{\alpha\sqrt{2n-1}\left(\frac{2n-1}{2n}\right)^{n-1}+\frac{1}{2^n-1}}<\alpha.\]
This is contradictory to the assumption $R(T^*)=\alpha$.\newline

\noindent $\bullet$~Case B ($0\le t_{N_*-1}^+$): Since $R$ is decreasing on $[0,t_{N_*-1}^+]$ and increasing on $[t_{N_*-1}^+,t_{N_*}^-]$,
\begin{align*}
	f(\alpha)-f(R(0))&=f(R(T^*))-f(R(0))\le f(R(T^*))-f(R(t_{N_*-1}^+)=\int_{A(t_{N_*-1})}^{A(T^*)}\frac{1}{\sin R}\frac{dR}{dA}dA\\
	&=\int_{2(N_*-1)+\frac{\pi}{2}}^{A(T^*)}\frac{1}{\sin R}\frac{dR}{dA}dA\le\int_{2(N_*-1)+\frac{\pi}{2}}^{A(T^*)}L_1(A)dA<\int_{-\frac{3\pi}{2}}^{-\frac{\pi}{2}}L_1(A)dA,
\end{align*}
which deduces the same contradiction in Case A.\newline

From the both cases above, we obtained $T^*\ge t_{N_*}^-$. Since we suppose $T^*<\infty$, there exists $l\in\mathbb{N}$ such that
\[t_l^-\le T^*\le t_{l+1}^-.\]
Since $R(t)$ is decreasing on $[t_l^-,t_l^+]$, one has $T^*>t_l^+$, which deduces
\begin{align*}
	f(R(T^*))-f(R(t_l^-))&=\big[f(R(t_l^+))-f(R(t_l^-))\big]+\big[f(R(T^*))-f(R(t_l^-))\big]\\
	&\le\int_{A(t_l^-)}^{A(t_l^+)}L_2(A)dA+\int_{A(t_l^+)}^{A(T^*)}L_1(A)dA\\
	&\le\int_{-\frac{\pi}{2}}^{\frac{\pi}{2}}L_2(A)dA+\int_{\frac{\pi}{2}}^{\frac{3\pi}{2}}L_1(A)dA<0.
\end{align*}
From this, we derive the contradiction
\[R(T^*)<R(t_l^-)<\alpha.\]
Therefore,
\[R(t)<\alpha,\quad\forall t\ge0.\]
\indent Now, we are ready to complete the proof using the above lemmas.
By Lemma \ref{L5.5} (ii),
\[R(t)\le R(t_l^-),\quad\mbox{for}~t\in[t_l^-,t_{l+1}^-].\]
Hence, a new step function $g$ which is defined as
\[g(t):=R(t_l^-)\quad t\in(t_l^-,t_{l+1}^-],\]
satisfies
\[R(t)\le g(t),\quad\forall t\ge0.\]
Also, one can derive
\[g(t)\le g(t_l^-),\quad\forall t\ge t_l^-.\]
Therefore, for $t\ge t_{N_*+1}^-$,
\begin{align*}
	\ln\frac{R(t)}{R(t_{N_*}^-)}&\le\ln\frac{g(t)}{R(t_{N_*}^-)}\le\ln\frac{R(t_{N_*+1}^-)}{R(t_{N_*}^-)}\\
	&\le-\frac{\sin\alpha}{\alpha}\cdot\frac{\Omega-2K}{2\pi}|C_2^+-C_1^-|\cdot|t_{N_*+1}^--t_{N_*}^-|\\
	&\le-\frac{\sin\alpha}{\alpha}\cdot\frac{\Omega-2K}{2\pi}|C_2^+-C_1^-|\cdot(t-t_{N_*}^-),
\end{align*}
which implies the desired result.
\end{proof}
\section{Numerical simulations} \label{sec:6}
\setcounter{equation}{0}
In this section, we provide several numerical simulations and compare them with analytical results in previous sections. For numerical simulations, we use the first-order forward Euler scheme with time-step $\Delta t = 10^{-2}$. 
\subsection{Temporal evolution of phase dynamics} \label{sec:6.1}
In this subsection, we performed several simulations for the phase dynamics of a single oscillator and coupled oscillators for higher-order couplings to observe dynamical patterns from different order $n$. For comparison, we used three different orders and the same other system parameters:
\[ n= 1, 10, 30, \quad (\nu, N, \kappa) = (5, 10, 1). \]
\subsubsection{Decoupled oscillator}  \label{sec:6.1.1} 
\begin{figure}
		\centering
		\begin{subfigure}[ht]{0.4\textwidth}
\centering
			\includegraphics[width=1.15\textwidth]{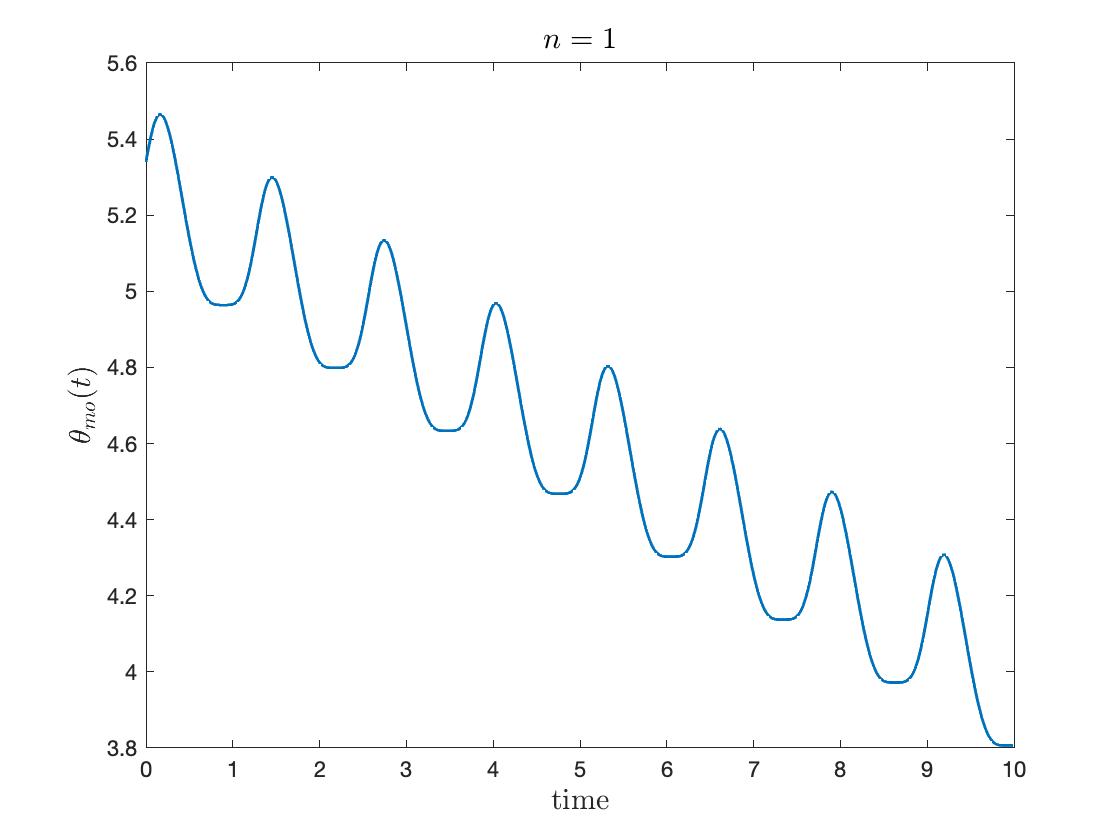}\caption{Dynamics of $\theta_{mo}$ for $n=1$} %\label{Fig-2-A}
		\end{subfigure}
\quad
		\begin{subfigure}[ht]{0.4\textwidth}
\centering
			\includegraphics[width=1.15\textwidth]{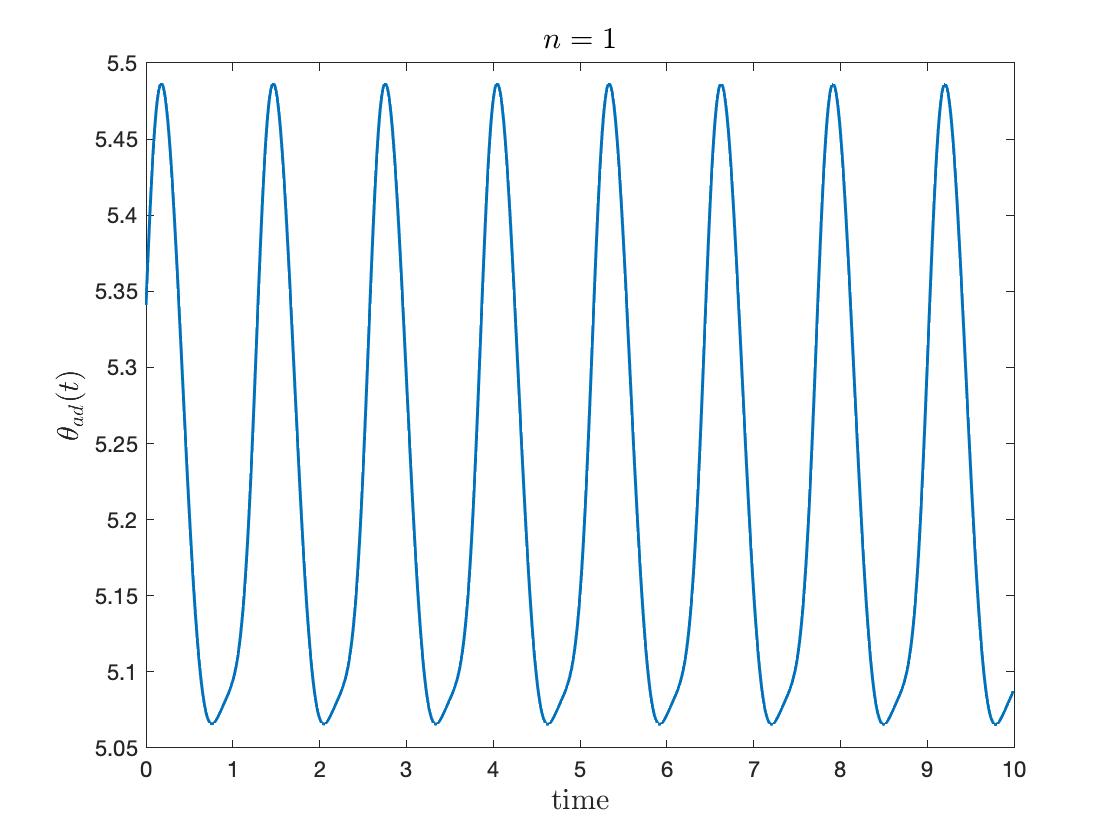}\caption{$\theta_{ad}$ with $\omega_{\infty}=0.1281$ and $n=1$} %\label{Fig-2-B}	
		\end{subfigure}
\\
\centering
		\begin{subfigure}[ht]{0.4\textwidth}
\centering
			\includegraphics[width=1.15\textwidth]{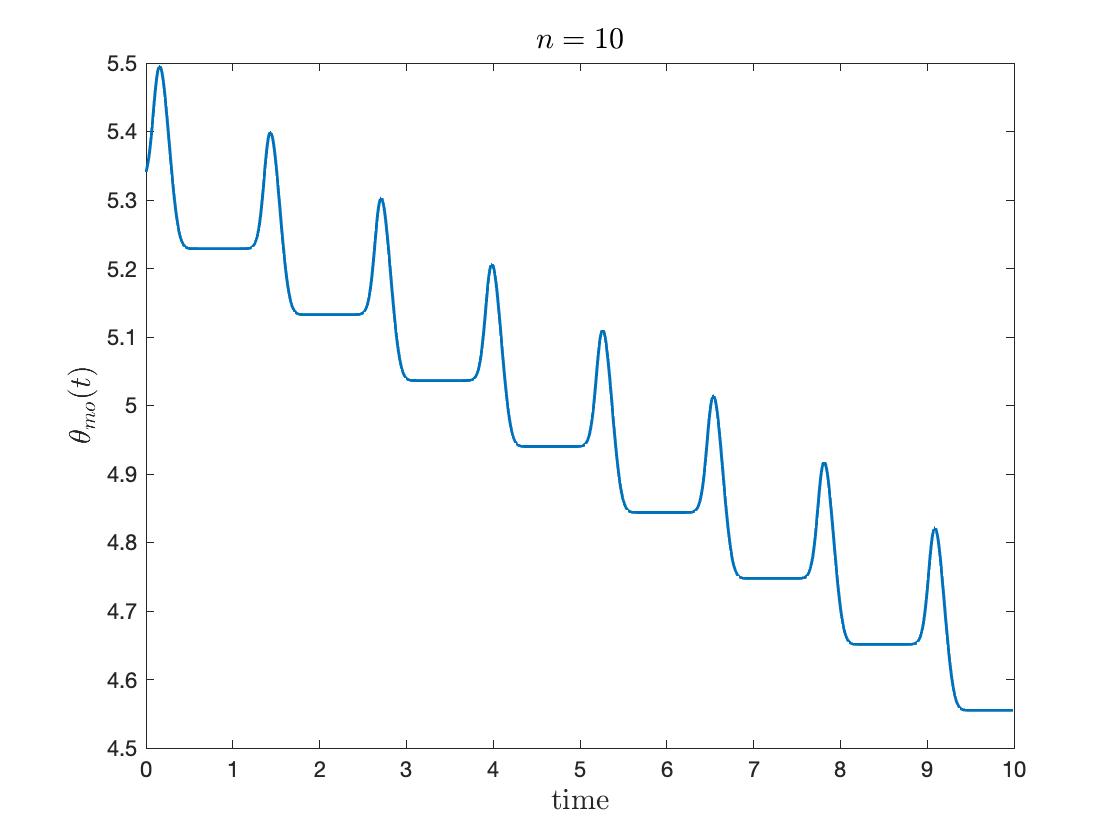}\caption{Dynamics of $\theta_{mo}$ for $n=10$}\label{Fig-3-A}
		\end{subfigure}
\quad
		\begin{subfigure}[ht]{0.4\textwidth}
\centering
			\includegraphics[width=1.15\textwidth]{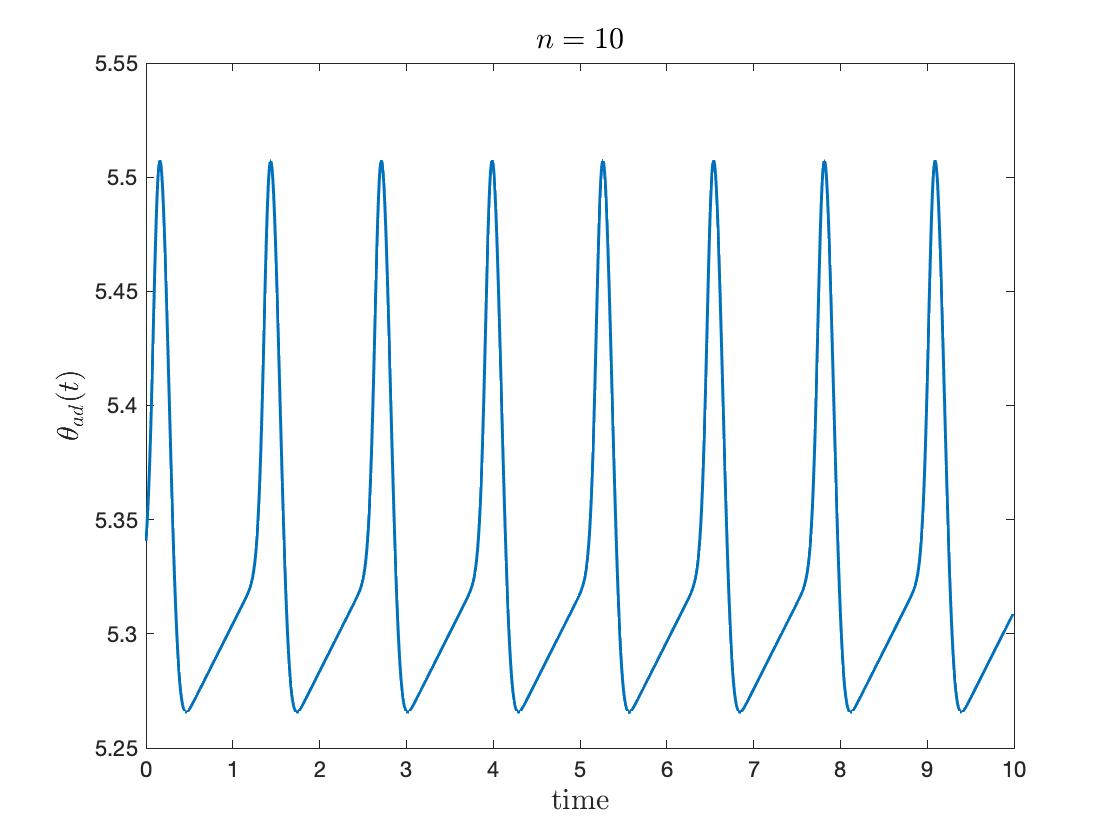}\caption{$\theta_{ad}$ with $\omega_{\infty}=0.0075$ and $n=10$}\label{Fig-3-B}	
		\end{subfigure}
\\
		\centering
		\begin{subfigure}[ht]{0.4\textwidth}
\centering
			\includegraphics[width=1.15\textwidth]{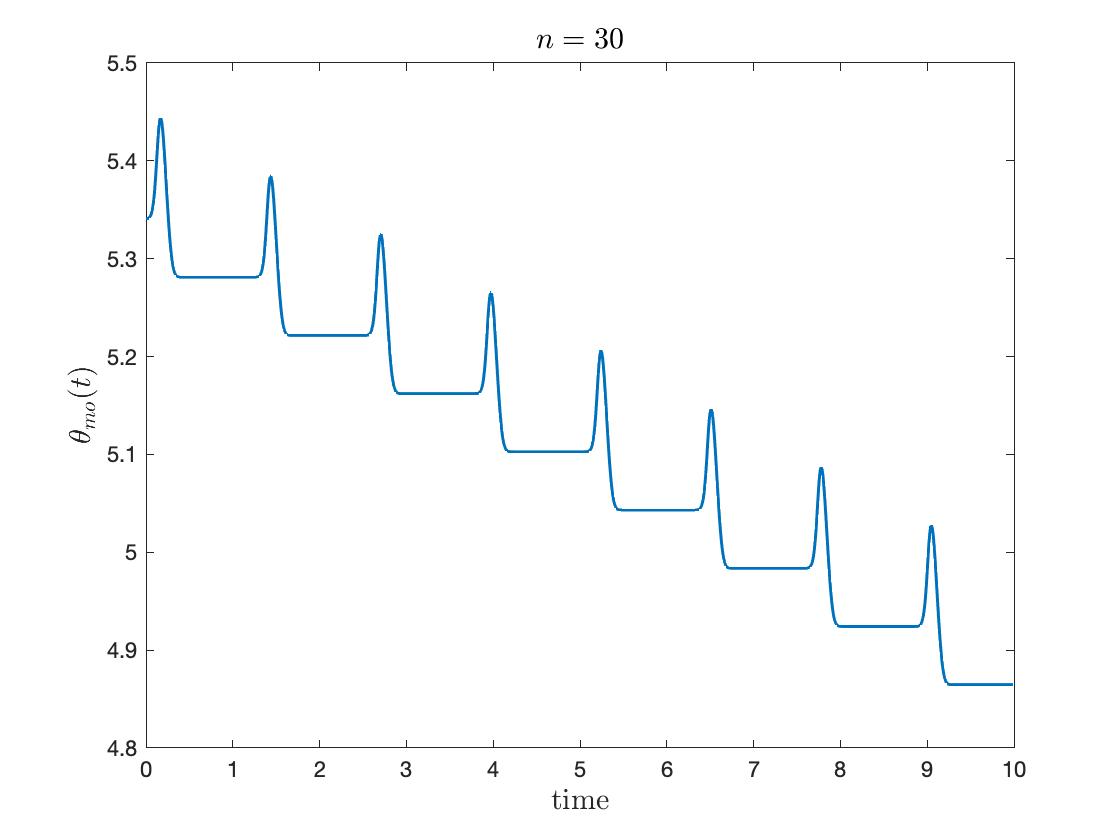}\caption{Dynamics of $\theta_{mo}$ for $n=30$}\label{Fig-4-A}
		\end{subfigure}
\quad
		\begin{subfigure}[ht]{0.4\textwidth}
\centering
			\includegraphics[width=1.15\textwidth]{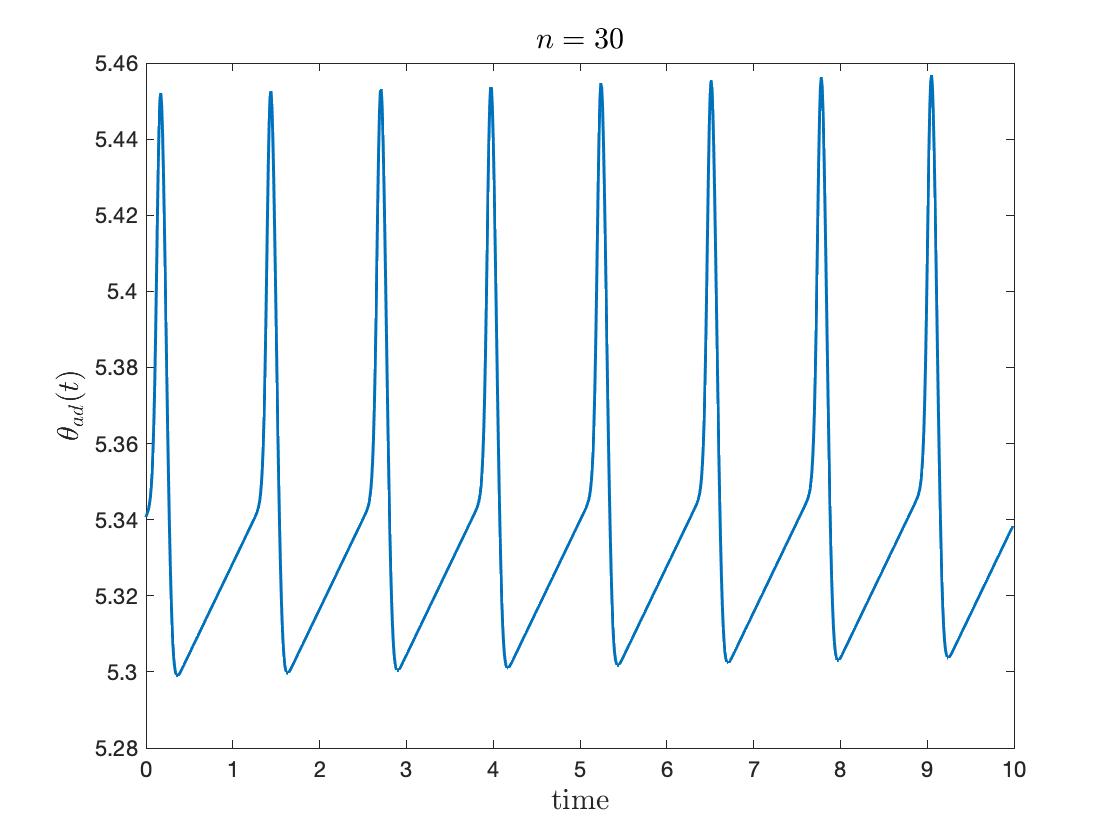}\caption{$\theta_{ad}$ with $\omega_{\infty}=0.0474$ and $n=30$}\label{Fig-4-B}	
		\end{subfigure}
		\caption{Single-oscillator dynamics}
		\label{Fig-4}
	\end{figure}

We first consider the dynamics of single oscillator with $\nu = 5$:
\begin{align*}
\begin{cases}
\displaystyle \dot{\theta} = 5 -  \Big( a_n  \left(1+\cos\theta\right)^n \Big) \sin \theta, \\
\displaystyle a_n=\frac{2n(2n-2)\cdots2}{2^n (2n-1)\cdots1}.
\end{cases}
\end{align*}
In Figure \ref{Fig-4}, we plotted temporal evolution of modified phases $\theta_{mo}$ and $\theta_{ad}$ with a suitable frequency shift:
\begin{align*}
\theta_{mo}(t):= \theta(t) - 5t, \quad   \theta_{ad}(t) := \theta(t) - (5 + \omega_{\infty}) t.
\end{align*}
Here $\omega_{\infty}$ is chosen after the simulation to nullify the drift from self-coupling. 
The spike-like phenomena are emphasized for large $n$, where the spike happens near $\theta(t) = 0$. 
\subsubsection{Coupled oscillators} \label{sec:6.1.2} 
\begin{figure}
		\centering
		\begin{subfigure}[ht]{0.4\textwidth}
\centering
			\includegraphics[width=1.15\textwidth]{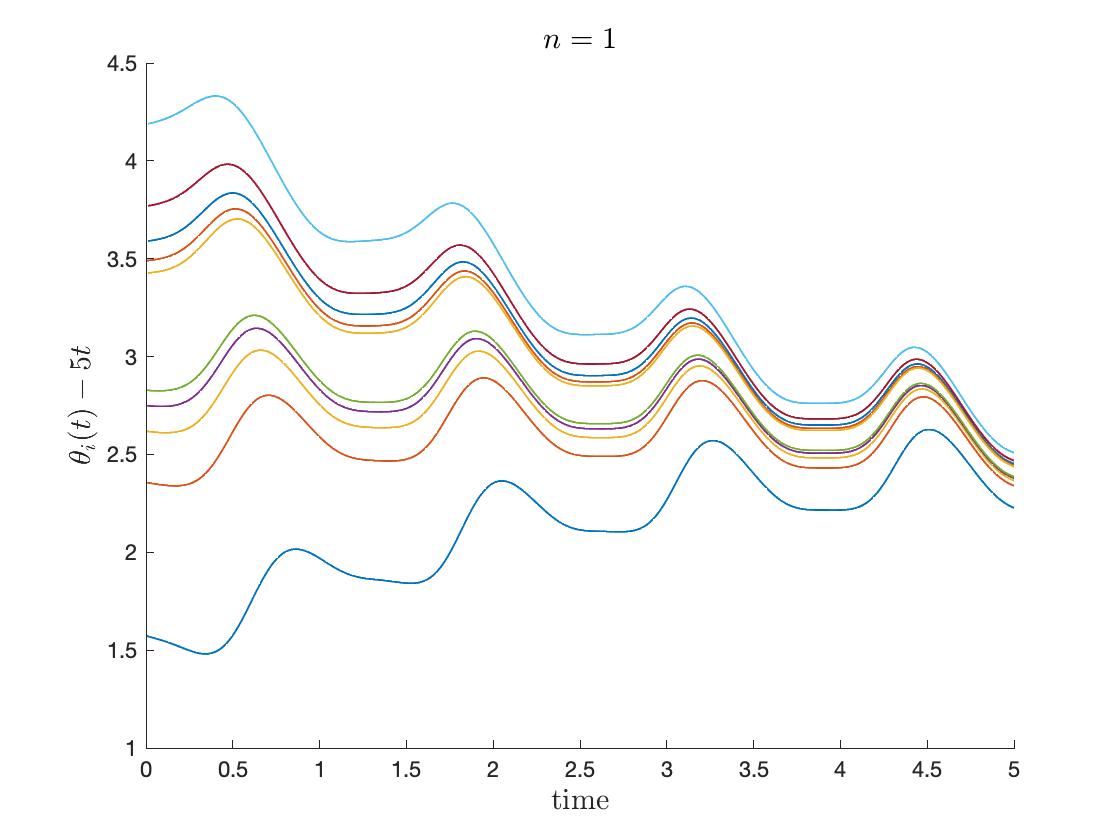}\caption{Dynamics of $\theta_i - 5t$ for $n=1$}\label{n1}
		\end{subfigure}
\quad
		\begin{subfigure}[ht]{0.4\textwidth}
\centering
			\includegraphics[width=1.15\textwidth]{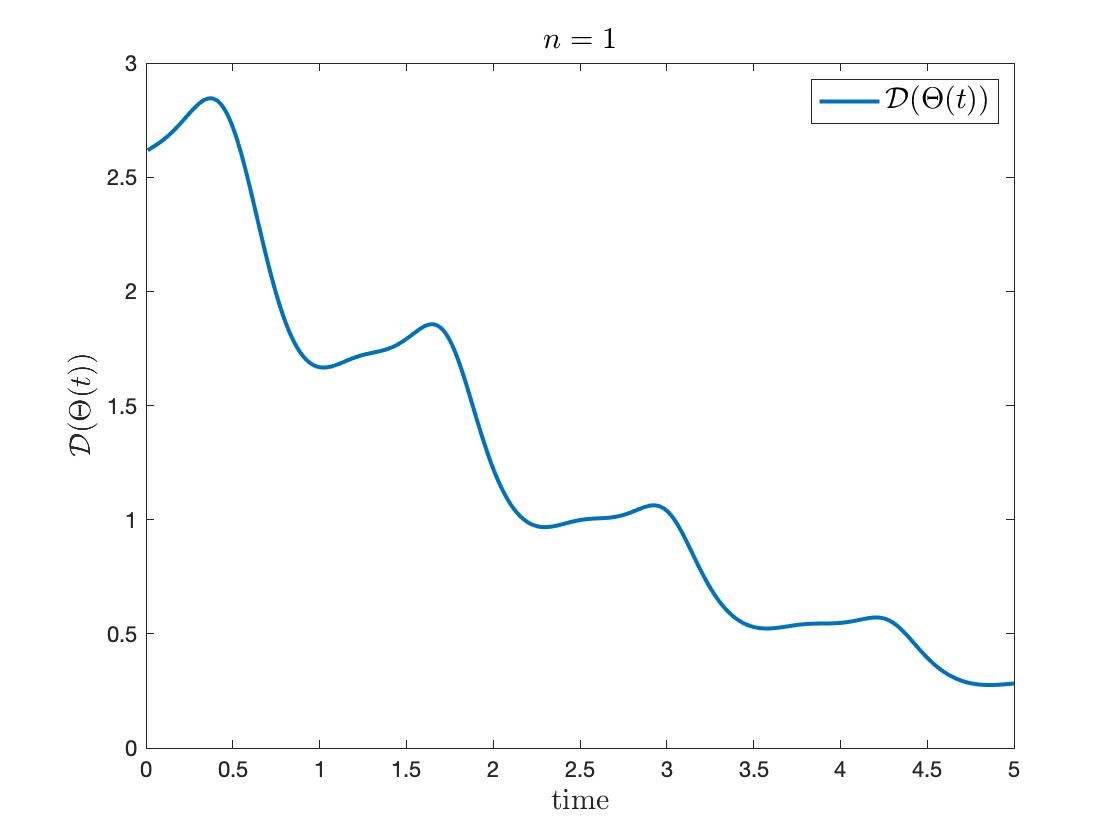}\caption{Dynamics of ${\mathcal D}(\Theta)$ for $n=1$}\label{n1D}
		\end{subfigure}
\vspace{.5cm}
		\centering
		\begin{subfigure}[ht]{0.4\textwidth}
\centering
			\includegraphics[width=1.15\textwidth]{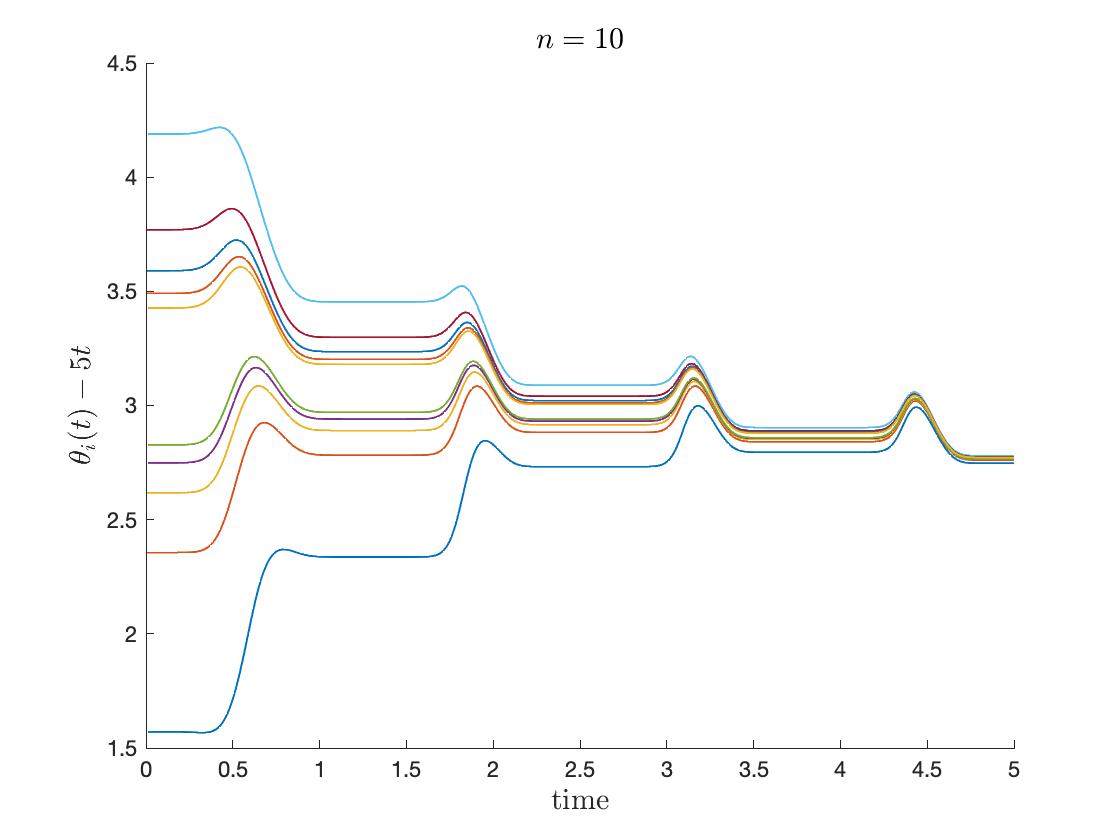}\caption{Dynamics of $\theta_i - 5t$ for $n=10$}\label{n10}
		\end{subfigure}
\quad
		\begin{subfigure}[ht]{0.4\textwidth}
\centering
			\includegraphics[width=1.15\textwidth]{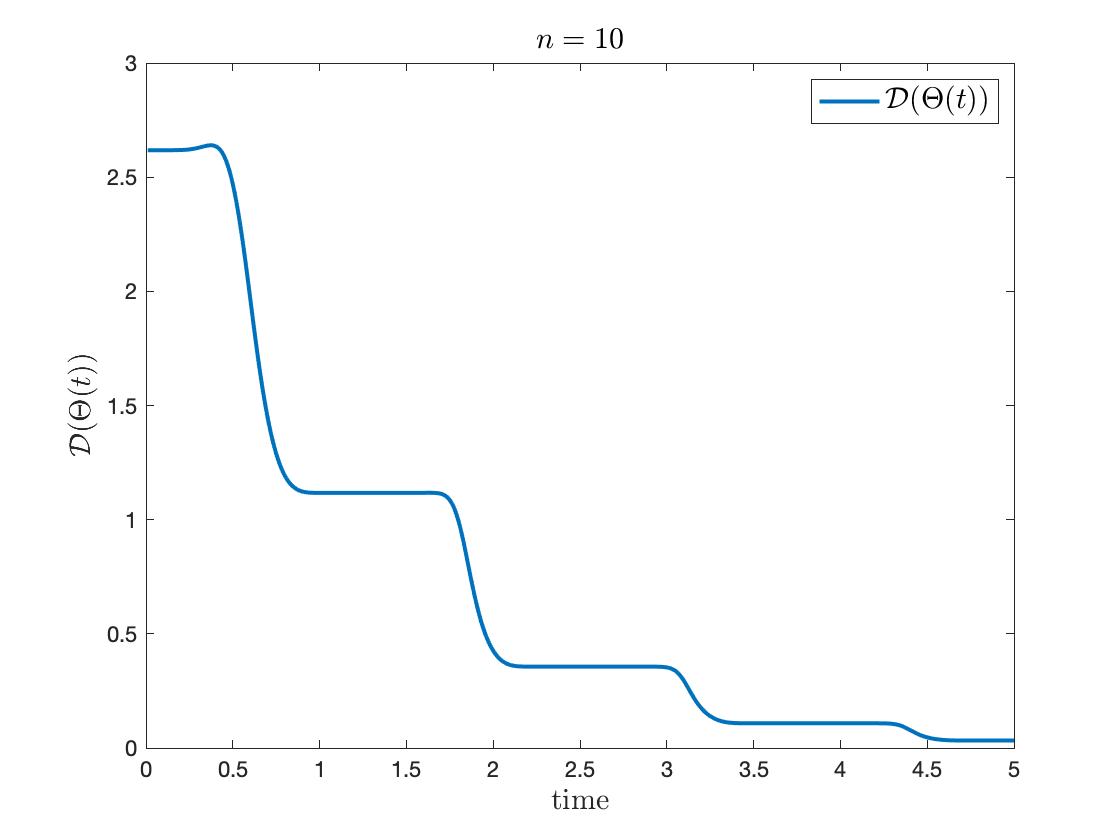}\caption{Dynamics of ${\mathcal D}(\Theta)$ for $n=10$}\label{n10D}
		\end{subfigure}
\vspace{.5cm}
		\centering
		\begin{subfigure}[ht]{0.4\textwidth}
\centering
			\includegraphics[width=1.15\textwidth]{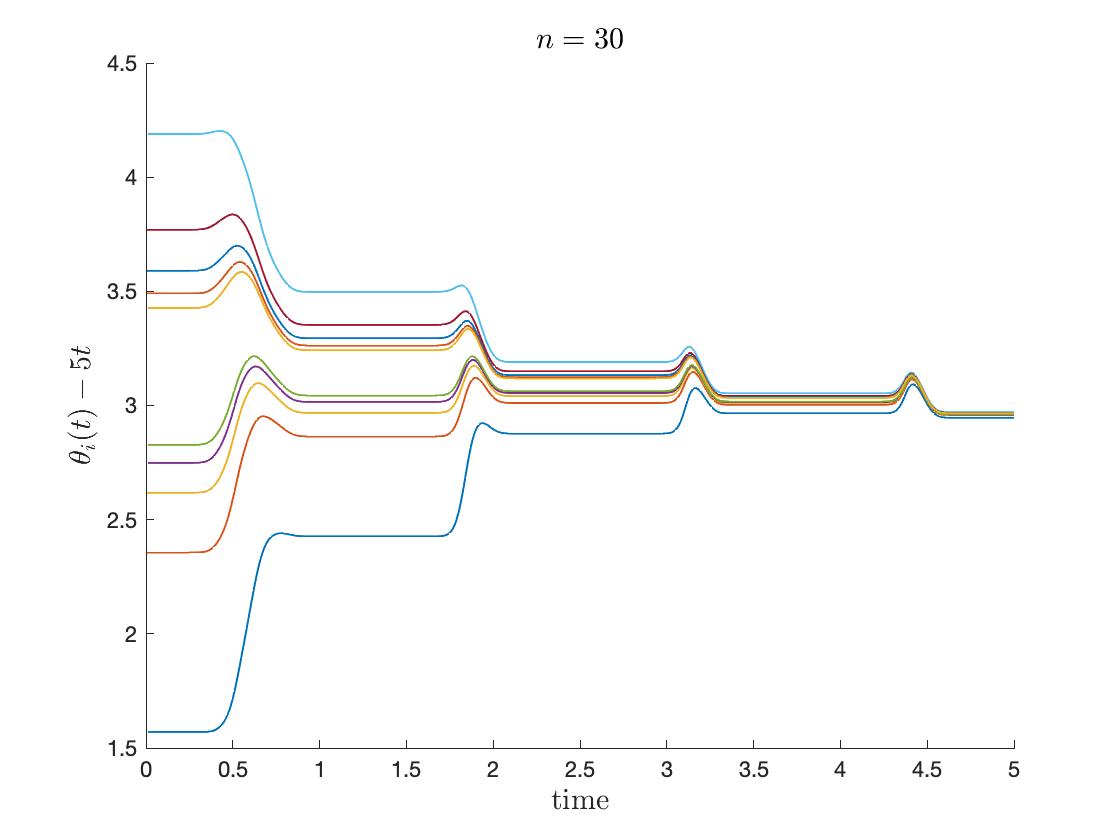}\caption{Dynamics of $\theta_i - 5t$ for $n=30$}\label{n30}
		\end{subfigure}
\quad
		\begin{subfigure}[ht]{0.4\textwidth}
\centering
			\includegraphics[width=1.15\textwidth]{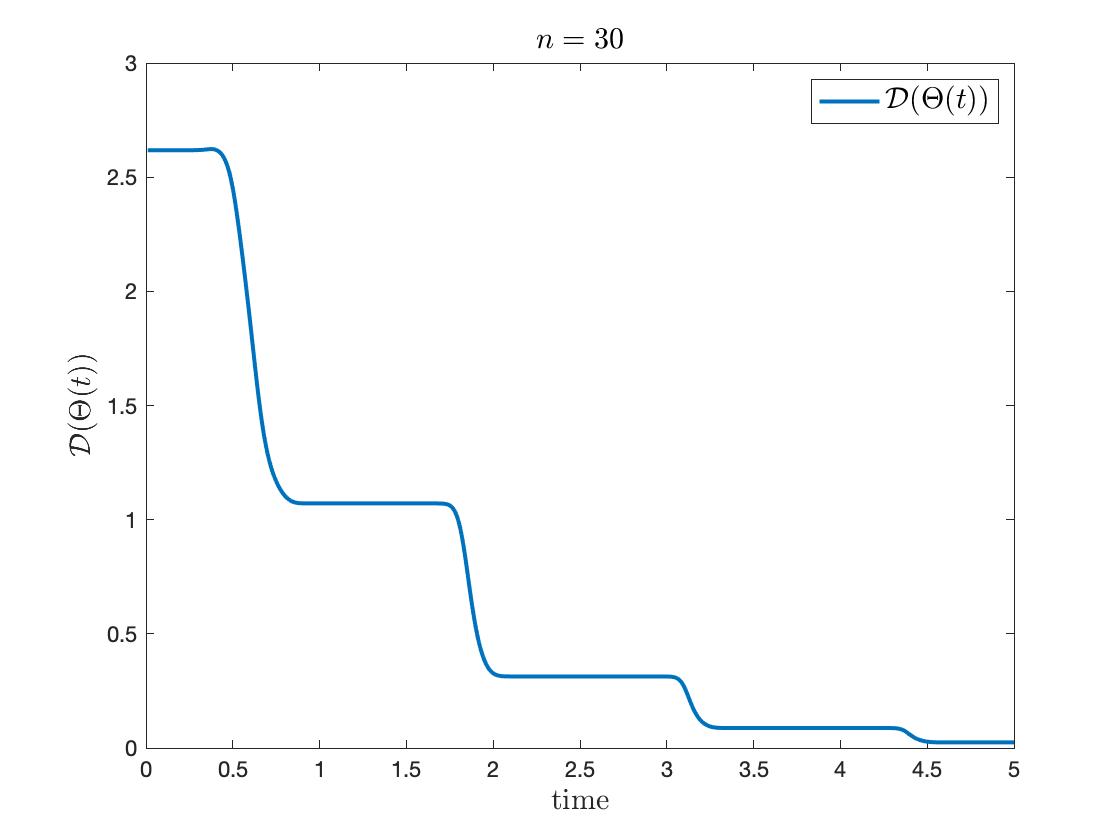}\caption{Dynamics of ${\mathcal D}(\Theta)$ for $n=30$}\label{n30D}	
		\end{subfigure}
		
		\caption{Coupled phase dynamics}
		\label{Fig-5}
	\end{figure}

Next, we consider the coupled dynamics for ten oscillators ($N=10$):
\begin{equation} \label{F-3}
\begin{cases}
\displaystyle \dot{\theta}_i= 5 -  \frac{1}{10} \Big( a_n \sum_{j=1}^{10} (1+\cos\theta_j)^n \Big) \sin\theta_i, \quad i \in [10], \\
\displaystyle a_n=\frac{2n(2n-2)\cdots2}{2^n (2n-1)\cdots1}.
\end{cases}
\end{equation}
\noindent In Figure \ref{Fig-5}, we use fixed initial data for different $n$:
\[ \{\theta_i^0\}_{i=1}^{10} = \Big \{ \frac{\pi}{2},~\frac{3\pi}{4},~\frac{5\pi}{6},~\frac{7\pi}{8},~\frac{9\pi}{10},~ \frac{4\pi}{3},~\frac{6\pi}{5},~\frac{8\pi}{7},~\frac{10\pi}{9},~\frac{12\pi}{11} \Big \}.
\]
We plotted the dynamics of modified phases $\theta_i - 5 t$ and phase diameter ${\mathcal D}(\Theta)$ over time, where all the numerical simulations tend to the locking state. %However, if you zoom in, you can see that phases oscillate. This increase is mainly due to the positivity of the natural frequencies. 
To observe the effect of interaction terms in \eqref{F-3}, we plot the following modified phases as in Figure \ref{Fig-4} by subtracting the drift from natural frequency $\nu = 5$:
	\[\theta_i(t)- 5 t,\quad i \in [10]. \]

From the whole simulations, there exist periodic oscillations. In detail, the motion of modified phases gets more accumulated near the spikes (where the oscillators `fire') as $n$ increases. For $n=1$, nearly the whole time the modified phases oscillate, while the diameter ${\mathcal D}(\Theta)$ shows stair-like plunges when $n=10$ or $30$. It suggests that the model acts more like a pulse-coupled dynamics as $n$ gets larger.

\subsection{Critical coupling strength} \label{sec:6.2}
\setcounter{equation}{0}
%	\begin{figure}
%		\centering
%		\begin{subfigure}[ht]{0.45\textwidth}
%			\centering
%			\includegraphics[width=1.2\textwidth]{gen_n1.jpg}\caption{Phases w/o natural frequencies}\label{gen_n1}
%		\end{subfigure}
%		\quad
%		\begin{subfigure}[ht]{0.45\textwidth}
%			\centering
%			\includegraphics[width=1.2\textwidth]{gen_n1D.jpg}\caption{Diameter of phases}\label{gen_n1D}	
%		\end{subfigure}
%		\caption{$n=1$}
%		\label{nis1_gen}
%	\end{figure}
%	\begin{figure}
%		\centering
%		\begin{subfigure}[ht]{0.45\textwidth}
%			\centering
%			\includegraphics[width=1.2\textwidth]{gen_n5.jpg}\caption{Phases w/o natural frequencies}\label{gen_n5}
%		\end{subfigure}
%		\quad
%		\begin{subfigure}[ht]{0.45\textwidth}
%			\centering
%			\includegraphics[width=1.2\textwidth]{gen_n5D.jpg}\caption{Diameter of phases}\label{gen_n5D}
%		\end{subfigure}
%		\caption{$n=5$}
%		\label{nis5_gen}
%	\end{figure}
%	\begin{figure}
%		\centering
%		\begin{subfigure}[ht]{0.45\textwidth}
%			\centering
%			\includegraphics[width=1.2\textwidth]{gen_n10.jpg}\caption{Phases w/o natural frequencies}\label{gen_n10}
%		\end{subfigure}
%		\quad
%		\begin{subfigure}[ht]{0.45\textwidth}
%			\centering
%			\includegraphics[width=1.2\textwidth]{gen_n10D.jpg}\caption{Diameter of phases}\label{gen_n10D}	
%		\end{subfigure}
%		\caption{$n=10$}
%		\label{nis10_gen}
%	\end{figure}
In Sections \ref{sec:3} to \ref{sec:5}, we have provided sufficient conditions on the coupling strength for incoherence, partial locking, locking and death. In the sequel, we numerically study missing ranges of coupling strength between these collective behaviors. 

More precisely, we study {\it ``empirical critical coupling strengths"} $\kappa_d,\kappa_p$ and $\kappa_i$ which correspond to the threshold values of coupling strength from death to phase-locking, from phase-locking to partial locking and from partial-locking to incoherence, respectively. For numerical values for $\kappa_d$, $\kappa_p$ and $\kappa_i$ depending on $n$, we choose the fixed initial data with $N=10$ for each scenario.

Firstly, we consider two scenarios with the nonidentical oscillators which corresponds to the situation in Section \ref{sec:3} and Section \ref{sec:4}. Here we use uniformly distributed natural frequencies
\[\{\nu_1,\nu_2,\cdots,\nu_{10}\}=\{5.4,5.8,\cdots,9.0\}\]
to distinguish the partial locking and locking significantly. The initial data are (uniformly) randomly chosen, one in a half circle and the other in the whole unit circle. 

\begin{figure}
		\centering
		\begin{subfigure}[ht]{0.4\textwidth}
			\centering
			\includegraphics[width=1.15\textwidth]{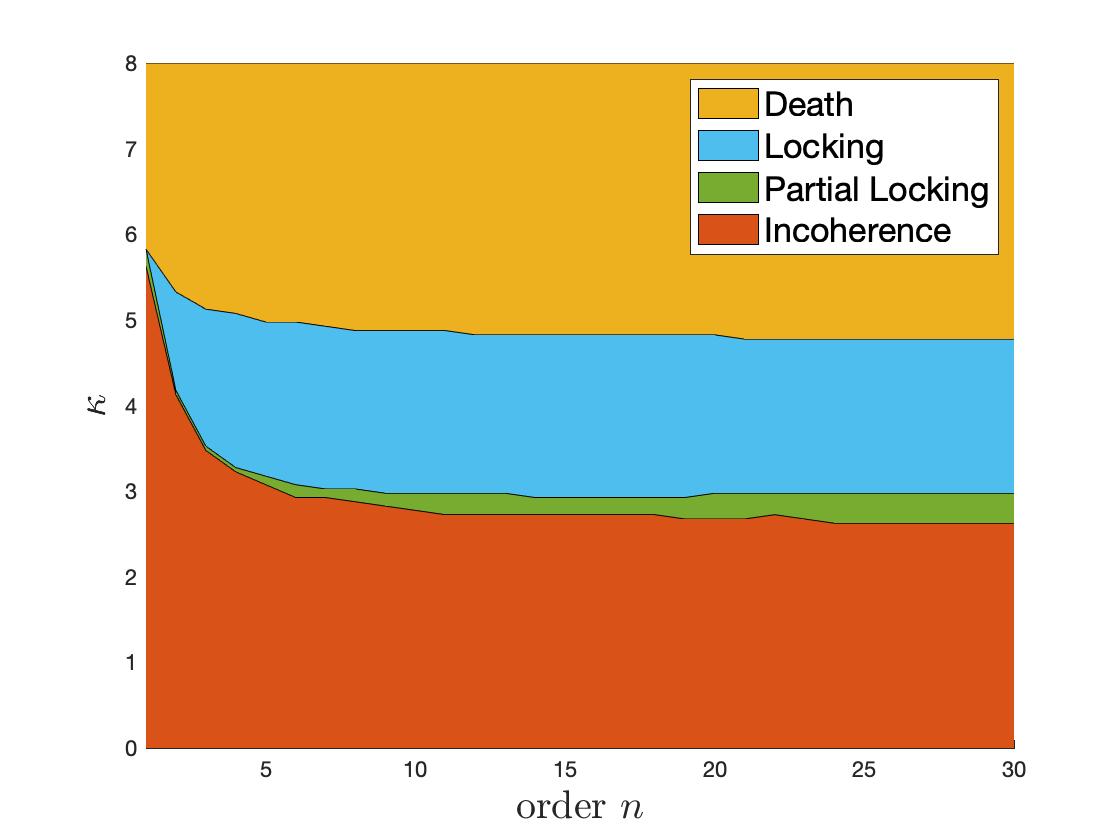}
			\caption{$\Theta(0)\in B\left(\frac{\pi}{2}\right)$}\label{nonid_1}
		\end{subfigure}
\quad
		\begin{subfigure}[ht]{0.4\textwidth}
			\centering
			\includegraphics[width=1.15\textwidth]{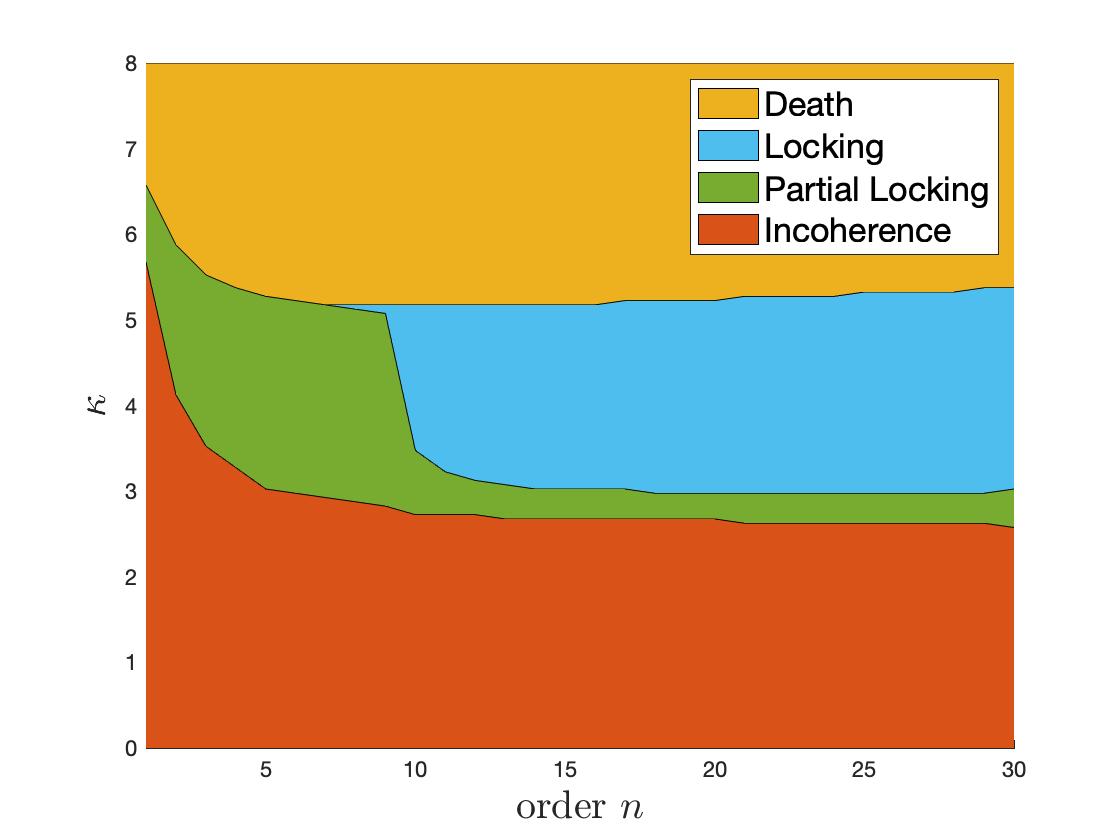}
			\caption{$\Theta(0)\in B\left(\frac{\pi}{200}\right)$}\label{nonid_2}
		\end{subfigure}
		\caption{Phase diagram with nonidentical oscillators}
		\label{Nonid}
\end{figure}

Figure \ref{Nonid} shows the phase diagram with different coupling strength $\kappa$ and order $n$. The asymptotic behaviors are labeled by observing the state at $t=500$ (with $\Delta t=0.01$) and marked with different colors. These simulations are done for each $n$ from $n=1$ to $n=30$ and each $\kappa$ from $0$ to $8$ with step size $\Delta\kappa=0.05$. 

One can observe that the asymptotic states depend on the initial data and $\kappa_d$ can not be said {\it increasing} or {\it decreasing} easily. We cannot directly apply the result of Section \ref{sec:4} since $\kappa_{d,n}$ is too large to compare to the value in the numeric phase diagram. However, the tendency along $n$ seems similar in the diagram while we expect $\kappa_{d} \sim {\mathcal O}(1)$. On the other hand, $\kappa_i$ decays with the rate next to $\mathcal{O}\left(\frac{1}{\sqrt{n}}\right)$, whose order was given in Section \ref{sec:3}. Overall, the numerical simulations shows similar tendency with respect to order $n$ to analytical results in Sections \ref{sec:3} and \ref{sec:4}. 

\begin{figure}
		\centering
		\begin{subfigure}[ht]{0.4\textwidth}
			\centering
			\includegraphics[width=1.2\textwidth]{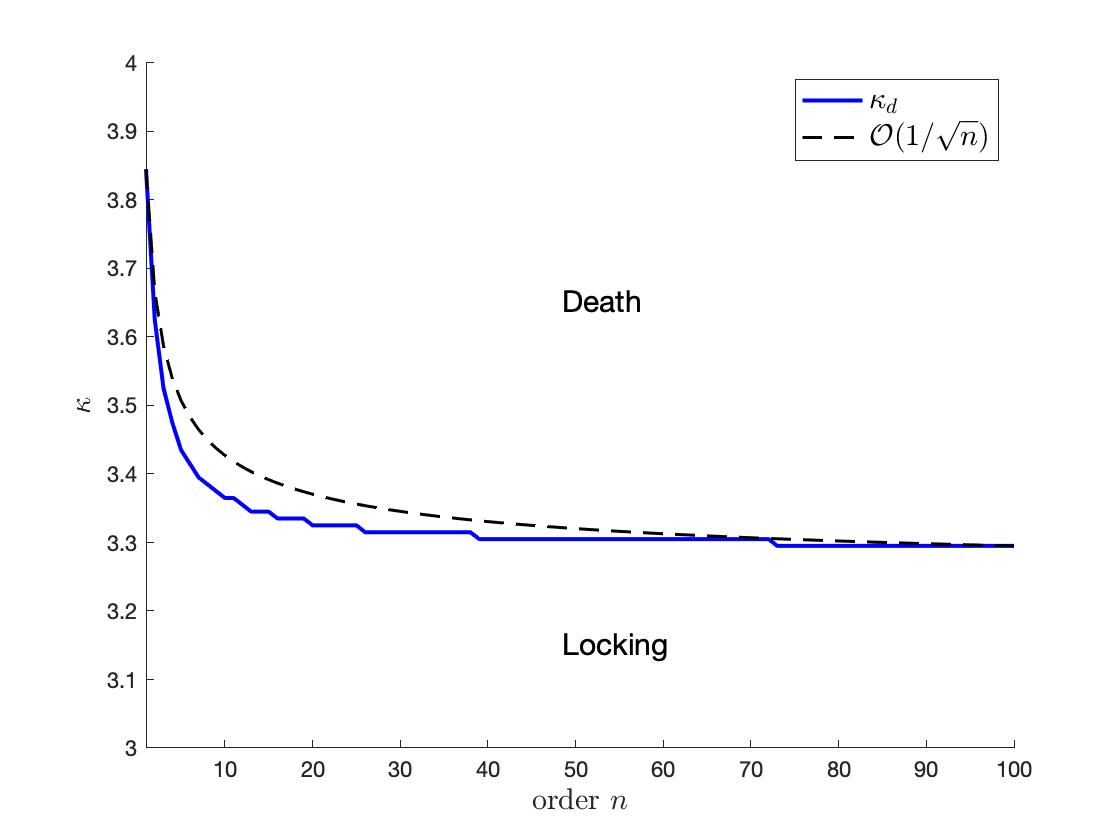}
			\caption{$\Theta(0)\in B\left(\frac{\pi}{2}\right)$}\label{id_half}
		\end{subfigure}
\quad
		\begin{subfigure}[ht]{0.4\textwidth}
			\centering
			\includegraphics[width=1.2\textwidth]{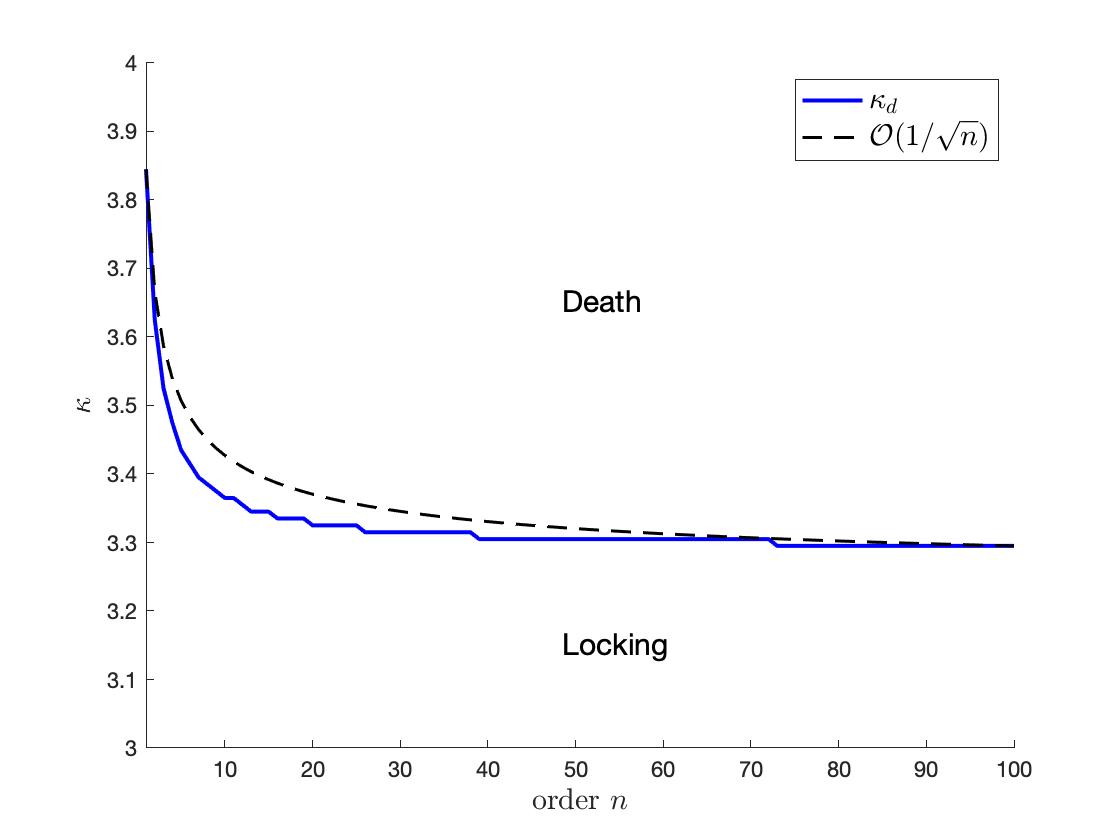}
			\caption{$\Theta(0)<B\left(\frac{\pi}{200}\right)$}\label{id_condition}
		\end{subfigure}
		\caption{Identical Oscillator}
		\label{Iden}
\end{figure}

Secondly, we also consider two scenarios with the identical oscillators (See Figure \ref{Iden}), $\nu=5$, corresponding to Section \ref{sec:5} using following settings:
\[t\in[0,500]~~\mbox{with}~~\Delta t=0.01,\quad n=1,\cdots,100,\quad\Delta\kappa =0.01.\]
Figure \ref{Iden}(A) shows the critical coupling $\kappa_d$ for the initial data in a half circle. For the identical oscillators, the dynamics only show locking or death since the whole ensemble will follow the common frequency without help of coupling interactions. The other one, Figure \ref{Iden}(B), is for the initial data which satisfies the condition in Theorem \ref{T5.1} for all $n$ with $\alpha = 3\pi/64$, hence, $\Theta(0)\in B(\pi/200)$.
Since these cases start with accumulated phases, one can expect that collective behaviors are much easier to emerge. Figure \ref{Iden} shows that, however, the empirical critical couplings strengths do not much depend on the initial data. As our estimate in Theorem \ref{T5.1}, $\kappa_d$ shows the monotone decrement in $n$ with a quite similar decaying rate that we analyze, $\mathcal{O}\left(\frac{1}{\sqrt{n}}\right)$. 

Therefore, we expect that our sufficient conditions on the coupling strength have the right growth rate with respect to order $n$.

\section{Conclusion}\label{sec:7}
\setcounter{equation}{0}
In this paper, we have provided several sufficient frameworks for the emergent dynamics of the Winfree model with higher-order couplings.
%For the first-order interaction, it is well-known that the Winfree model exhibits asymptotic patterns such as incoherence, locking and death. Thus, it is reasonable to ask how the order of interactions affects aforementioned asymptotic patterns. The Winfree model with higher-order couplings can be understood as an approximate model for the pulse-coupled synchronization model. The Integrate-and-fire model is a prototype pulse-coupled model for synchronous dynamics of neurons, and it has been extensively used in neuroscience. However, 
Rigorous study on the pulse-coupled synchronous dynamics is very rare compared to the phase-coupled models such as the Kuramoto model and the Winfree model. 
As a first step toward the mathematical analysis of pulse-coupled models, we adopted the Winfree model with higher order trigonometric couplings, and study its emergent dynamics with respect to the order of couplings $n$. As the order $n$ increases, the influence function approaches to the constant multiple of Dirac mass. From the estimate on phase diameters, we analyzed how the order $n$ of influence function affects the emergence of various types of collective modes. Of course, our proposed frameworks are only sufficient ones, and they certainly not optimal at all. It would be interesting to derive more refined and optimal coupling strengths between the transitions between collective modes:
\begin{center}
Incoherence $\quad \Longrightarrow \quad$ locking $\quad \Longrightarrow \quad$ death.
\end{center}

On the other hand, compared to the numerical simulations in previous section, our coupling strength condition from Theorem \ref{T3.1} to Theorem \ref{T5.1} only work for restricted initial data, not for general one as in Figure \ref{Nonid}(A) or \ref{Iden}(A). From these simulations, one may expect that the same asymptotic behavior emerges under a similar range of coupling strength $\kappa$. We leave them as a future work.

\end{document}